\documentclass{aims_emi}

\usepackage{txfonts}
\def\typeofarticle{Preprint}
\def\currentvolume{}
\def\currentissue{}

\def\ppages{}
\def\DOI{}

\def\Revised{}
\def\Accepted{}
\def\Published{}

\usepackage{amsmath,amsfonts,amssymb,bbm,amsthm}
\usepackage{mathrsfs}
\usepackage{graphicx}  
\usepackage{geometry} 
\usepackage{overpic}
\usepackage{algorithm,algpseudocode} 
\usepackage{tikz}\usetikzlibrary{patterns}
\usepackage{xcolor}
\usepackage{bbold}

\newcommand{\R}{\mathbb R}
\newcommand{\Dx}{\Delta x}
\newcommand{\Dt}{\Delta t}
\newcommand{\demand}{\textsc{d}}
\newcommand{\supply}{\textsc{s}}
\newcommand{\mis}[1]{\mu^{#1}}
\newcommand{\miss}[1]{\mathring{\mu}^{#1}}
\newcommand{\xmin}{x_{\textsc{min}}}
\newcommand{\xmax}{x_{\textsc{max}}}
\newcommand{\rjh}{r^\supply_{jh}}
\newcommand{\skl}{r^\demand_{k\ell}}
\newcommand{\gammacompleto}{\gamma_{(jh,k\ell)}}
\newcommand{\rhouno}{\rho^\supply}
\newcommand{\rhodue}{\rho^\demand}
\newcommand{\muuno}{\mu^\supply}
\newcommand{\mudue}{\mu^\demand}
\newcommand{\Funo}{F^{\supply}}
\newcommand{\Fdue}{F^{\demand}}
\newcommand{\rhomax}{\rho_{\textsc{max}}}
\newcommand{\velmax}{v_{\textsc{max}}}
\newcommand{\rhoin}{\rho_{\textsc{in}}}
\newcommand{\rhoout}{\rho_{\textsc{out}}}
\newcommand{\vin}{v_{\textsc{in}}}
\newcommand{\vout}{v_{\textsc{out}}}
\newcommand{\rhomeno}{\rho_-}
\newcommand{\rhopiu}{\rho_+}
\newcommand{\fin}{F_{\textsc{in}}}
\newcommand{\fout}{F_{\textsc{out}}}
\newcommand{\pfin}{\hat F_{\textsc{in}}}
\newcommand{\pfout}{\hat F_{\textsc{out}}}
\newcommand{\mres}{\mathbin{\vrule height 1.6ex depth 0pt width
		0.13ex\vrule height 0.13ex depth 0pt width 1.3ex}}
\newcommand{\unoN}{\mathbb{1}_{\!N}}
\newcommand{\unoNq}{\mathbb{1}_{\!N^2}}
\newcommand{\zeroN}{\mathbb{0}_{\!N}}
\newcommand{\me}{\mathbf e}
\newcommand{\fg}{\textsc{fg}}
\newcommand{\sss}{\textsc{ss}}
\newcommand{\pr}{\textsc{pr}}
\newcommand{\ghk}{\textsc{ghk}}

\newtheorem{dfn}{Definition}
\newtheorem{rmk}{Remark}

\numberwithin{equation}{section}

\makeatletter
\renewcommand{\@biblabel}[1]{#1\hfill \hspace{-0.2cm}}
\makeatother

\begin{document}

\title{Numerical Computation of Generalized Wasserstein Distances with Applications to Traffic Model Analysis}

\author{%
  Maya Briani\affil{1},
  Emiliano Cristiani\affil{1,}\corrauth,
  Giovanni Franzina\affil{1}
  and
  Francesca L.\ Ignoto\affil{1,2,}\corrauth
}

\shortauthors{the Author(s)}

\address{%
  \addr{\affilnum{1}}{Istituto per le Applicazioni del Calcolo, Consiglio Nazionale delle Ricerche, Rome, Italy.}
  \addr{\affilnum{2}}{Dipartimento di Scienze di Base e Applicate per l'Ingegneria, Sapienza Universit\`a di Roma, Rome, Italy.}
}

\corraddr{e.cristiani@iac.cnr.it; francescalourdes.ignoto@uniroma1.it}

\begin{abstract}
Generalized Wasserstein distances allow to quantitatively compare two continuous or atomic mass distributions with equal or different total mass. 
In this paper, we propose four numerical methods for the approximation of three different generalized Wasserstein distances introduced in the last years, giving some insights into their physical meaning. 
After that, we explore their usage in the context of the sensitivity analysis of differential models for traffic flow. 
The quantification of models sensitivity is obtained by computing the generalized Wasserstein distances between two (numerical) solutions corresponding to different inputs, including different boundary conditions. 
\end{abstract}


\maketitle
\tableofcontents
\section{Introduction}\label{sec:intro}

\paragraph{Context.} 
In this paper, we are concerned with the numerical computation of Wasserstein distance (WD) and Generalized Wasserstein distance (GWD) in the context of traffic flow models.
As for all dynamical systems, the study of the sensitivity of the models can be realized measuring the `distance' between two solutions obtained with different inputs (like, e.g., any of the model parameters, initial conditions, boundary conditions, etc.). 
This allows one to understand their impact on the final solution, and ultimately quantify the degree of chaoticity of the system. 
The question arises which distance is more suitable to this kind of investigation. 
It is by now well understood that $L^p$ distances do not catch the natural concept of distance among traffic (vehicle) densities, see, e.g., the discussion in \cite[Section 7.1]{cristianibook} and Section \ref{sec:WD} in this paper, while the WD appears more natural, at least in the context of traffic flow, see, e.g., 
\cite{
	briani2018CMS, 
	cristianibook, 
	cristiani2019M2AS}.
The drawback is that WD is limited to \emph{balanced} mass distribution (equal masses), while real traffic problems often need to consider scenarios with a different amount of vehicles, especially because of different inflow/outflow at boundaries. 
This suggests to move towards GWDs, which allows to deal with \emph{unbalanced} mass distributions, but serious computational problems arise: all existing definitions are quite abstract and it is not immediate to derive suitable numerical approximations.

\paragraph{Relevant literature.}
A common feature of various optimal transport
distances between finite positive Borel measures
on a separable complete metric space $(X,\mathfrak{D})$
is that
their primal definition requires the minimization of a
transportation cost including an additive contribution of the  form
\begin{equation}
	\label{eq:cost}
	\int_{X\times X} c(x,y)\,d\gamma\,,
\end{equation}
where $c:X\times X\to[0,+\infty]$ is a fixed lower semi-continuous cost function.

In the balanced problem, that one considers to define the classical Wasserstein distance between two probability measures $\mu^\supply$ and $\mu^\demand$, the minimization occurs upon transport plans $\gamma$, i.e., probability
measures on $X\times X$  with marginals $\mu^\supply$ and $\mu^\demand$. Namely, one
considers the projection $\pi^{\supply}$
(resp., $\pi^{\demand}$)
onto the first (resp., the second) factor of the Cartesian product $X\times X$, 
and takes as admissible competitors the finite positive
measures $\gamma$ for which 
\begin{equation}
	\label{eq:balanced}
	\pi^\supply_{\,\sharp}\gamma=\mu^\supply\,,\qquad
	\pi^\demand_{\,\sharp}\gamma=\mu^\demand\,,
\end{equation}
where $\pi^\supply_{\,\sharp}\gamma$ (resp.,
$\pi^\demand_{\,\sharp}\gamma$) denotes the 
push-forward measure of $\gamma$ under $\pi^{\supply}$
(resp., $\pi^\demand)$.

In order to define distances
between finite positive measures on a bounded open set $\Omega$
in a finite dimensional Euclidean space,
there exist various approaches to similar, but \emph{unbalanced}, problems, that
correspond each to the appropriately constrained minimization
of a specific objective functional.
In this paper we focus on four of such approaches.
\begin{enumerate}
	\item The approach by Figalli \& Gigli \cite{figalli2010JMPA}, which consists in minimizing the integral 
	cost \eqref{eq:cost}
	among all finite positive measures $\gamma$ on 
	$\overline\Omega\times\overline\Omega$ whose marginals
	induce on $\Omega$ the prescribed finite positive  Borel measures. 
	\item The approach by Piccoli \& Rossi~\cite{piccoli2014ARMA}, which includes in the definition
	of the functional the additive contributions
	\begin{equation}
		\label{eq:PR}
		|\pi^\supply_{\,\sharp}\gamma-\mu^\supply|(\Omega)
		+
		|\pi^\demand_{\,\sharp}\gamma-\mu^\demand|(\Omega)\,,
	\end{equation}
	in addition to \eqref{eq:cost}, with the minimization
	occurring upon finite positive Borel measures
	on $\Omega\times\Omega$.
	\item A third approach
	relies on
	the logarithmic entropy-transport
	distance
	that was independently introduced 
	by Kondratyev, Monsaingeon \& Vorotnikov in~\cite{KonMorVor2016AdvDE},
	by Chizat, Peyr\'e, Schmitzer \& Vialard in~\cite{ChiPeySchVia2018FoundCompMath,ChiPeySchVia2018JFA},
	and by Liero, Mielke \& Savar\'e in~\cite{LieMieSav2018InventMath},
	called Hellinger--Kantorovich or Wasserstein--Fisher--Rao distance. This
	approach also relies
	on the minimization, among finite positive Borel measures
	on $\Omega\times\Omega$, of a perturbation of \eqref{eq:cost},
	except that the penalization \eqref{eq:PR} is 
	replaced by
	\begin{equation}
		\label{eq:SS}
		\int_\Omega 
		\Big(
		\sigma^\supply \log(\sigma^\supply) -\sigma^\supply+1
		\Big)
		\,d\mu^\supply
		+
		\int_\Omega 
		\Big(
		\sigma^\demand \log(\sigma^\demand) -\sigma^\demand+1
		\Big)
		\,d\mu^\demand\,,
	\end{equation}
	if $d(\pi^\supply_{\,\sharp}\gamma)=\sigma^\supply\,d\mu^\supply$
	and $d(\pi^\demand_{\,\sharp}\gamma)=\sigma^\demand\,d\mu^\demand$
	for appropriate densities $\sigma^\supply$ and
	$\sigma^\demand$, and $+\infty$ otherwise.
	(Incidentally, a similar objective, in which an entropic regularization
	is added to \eqref{eq:cost} and \eqref{eq:SS},
	in considered in~\cite{ma2021AAAI}, where
	the appropriate dual formulation of the problem
	is addressed by Sinkhorn iterations.)
	\item A fourth approach
	relates to the contribution by Savar\'e \& Sodini~\cite{savare2024JMPA},
	in which the logarithmic entropy transport distance
	is characterized, in its various equivalent formulations, by convex relaxation. In particular, their approach can be considered
	the natural one to numerically approximate the (Gaussian) Hellinger--Kanrotovich distance in its homogeneous formulation.
\end{enumerate}
We point out that the last two approaches
only differ numerically, as in
principle the corresponding
unbalanced optimal transport formulations are 
analytically equivalent. Both extend WD to the cone of non-negative Borel measures of finite mass that are different from zero (see Remark \ref{GHKsiext} below). 
On the contrary,
the first two approaches do give an extension
of WD, despite defining a distance on the cone (see Remarks \ref{FGnotext} and \ref{PRnotext} below). 
We recall that no extension of WD to the set of signed measures, possibly zero, can define a distance~\cite{lombardini2022PAMS} (non-canonical optimal transport distance have however
been successfully defined, see~\cite{Amb-Mai-Ser2011,Mai2012}).

\paragraph{Main contribution.} 
This paper has multiple goals:
first of all, the four approaches are recalled and compared at the theoretical level, so to have a common and exhaustive reference for all recent developments about GWDs.

Second, four new \emph{ad hoc} computational methods are devised in order to practically compute the four GWDs (actually they are three GWDs, since two of them coincide under suitable assumptions).
The methods are individually studied and their main properties are discussed through some explanatory numerical tests.

Third, the four approaches are compared on the same academic test in order to understand which one is the most suitable in the context of traffic flow models.

Finally, the three most promising approaches are further compared on four realistic traffic problems.

Anticipating the conclusions, we can say that Figalli \& Gigli approach appears to be preferable among all, mainly for two reasons: (i) it is the most manageable from the computational point of view (easy to code, fast computation), and (ii) it is the most informative on the modeling side.

\paragraph{Plan of the paper.}
The paper is organized as follows. 
In Section \ref{sec:WD} we recall the main motivations of the paper and we set up notations for linear programming and discrete-to-discrete optimal mass transport.
In Section \ref{sec:GWD} we recall the definition and the basic properties of some known
unbalanced optimal transport-type distances between positive Borel measures of finite mass.
In Section \ref{sec:GWDcomputational} we revisit the distances 
in the discrete-to-discrete case and, more important, we propose computational approaches for their approximation. The distances are then compared on an academic numerical test. In addition -- postponed in the Appendix -- we present some tests specifically designed for each distance in which the exact value is available analytically, in order to study their interplay with the numerical boundary and their capability to detect significant patterns.
In Section \ref{sec:traffic} the proposed methods are applied in the context of traffic analysis, in order to achieve various goals: 
(i) to assess the impact of boundary conditions, 
(ii) to evaluate the sensitivity to traffic lights, 
(iii) to compare forecasts depending on whether or not the model takes inertial effects into account, and 
(iv) to compare the usage of the distances in calibration and validation phases. 
	The paper ends with a final discussion.

\section{The Wasserstein distance in the context of vehicular traffic}\label{sec:WD}

\subsection{Motivations}
\begin{figure}[h!]
	\centerline{
		\begin{overpic}
			[width=0.4\textwidth]{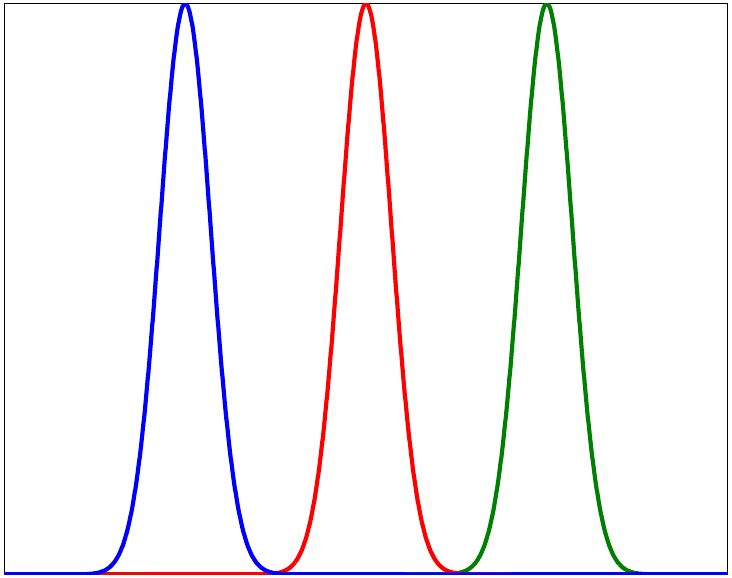}
			\put(17,60){$\rho^1$}\put(41.5,60){$\rho^2$}\put(66,60){$\rho^3$}
	\end{overpic}} 
	\caption{Three density functions with disjoint supports.}
	\label{fig:L1vsW}
\end{figure}

Usually the quantification of the closeness of two time-dependent distributions is performed by means of the $L^1$, $L^2$ or $L^\infty$ distance (in space at final time or in both space and time). Although this can be satisfactory for nearly equal outputs or for convergence results, it appears inadequate for measuring the distance of largely different outputs. 
To see this, let us focus on the example depicted in Figure \ref{fig:L1vsW}, where 
three density functions $\rho^i$, $i=1,2,3$, corresponding to the same total mass, say $M$, are plotted. It is plain that the $L^1$ distances between $\rho_1$ and $\rho_2$ and between $\rho_1$ and $\rho_3$ are both equal to $2M$. Similarly, all $L^p$ distances are blind with respect to variation of the densities once the supports of them are disjoint. Our perception of distance suggests instead that $\|\rho_3-\rho_1\|>\|\rho_2-\rho_1\|$, and this is exactly what Wasserstein distance guarantees, as we will see in the next section.

\subsection{Basic theory}\label{sec:WD-teoria}

\noindent 
Let us denote by $(X,\mathfrak{D})$ a locally compact complete and separable metric space with distance $\mathfrak{D}$, and by $\mathcal B(X)$ the Borel $\sigma$-algebra of $(X,\mathfrak{D})$. 
Let us also denote by $\mathcal M^+(X)$ the set of non-negative finite measures on $(X,\mathcal B(X))$. 
Let $\mis{\demand}$ ($\demand$ standing for demand or destination) and $\mis{\supply}$ ($\supply$ standing for supply or source) be two Radon measures in $\mathcal M^+(X)$ such that $\mis{\demand}(X)=\mis{\supply}(X)$ (same total mass).
\begin{dfn}[Wasserstein distance] 
	For any $p\in[1,+\infty)$, the $L^p$-Wasserstein distance between $\muuno$ and $\mudue$ is
	\begin{equation}\label{def:W}
		W_p(\muuno,\mudue):=\left(\inf_{\gamma\in\Gamma(\muuno,\mudue)}\int_{X\times X}\mathfrak{D}(x,y)^p\ d\gamma(x,y)\right)^{1/p}
	\end{equation}
	where $\Gamma$ is the set of transport plans connecting $\muuno$ to $\mudue$, i.e.
	\begin{equation*}
		\Gamma(\muuno,\mudue):=\left\{\gamma\in\mathcal Q^+(X\times X)\text{ s.t. }\gamma(A\times X)=\muuno(A),\ \gamma(X\times B)=\mudue(B),\ \forall \ A,B\in \mathcal{B}(X)\right\}.
	\end{equation*}
\end{dfn}
\begin{rmk}\label{rmk:Omega}
	The distance defined by \eqref{def:W}
	makes sense for any two
	given measures $\muuno$ and $\mudue$
	on a finite dimensional normed space,
	whether they are absolutely continuous
	or not with respect to the Lebesgue measure.
\end{rmk}

It is well known that
the notion of Wasserstein distance can be put in relation with the Monge--Kantorovich optimal mass transfer problem \cite{villani2009book}: a pile of, say, soil, has to be moved to an excavation with same total volume. Moving a unit
quantity of mass has a cost which equals the distance between the source and the
destination point. We consequently are looking for a way to rearrange the first mass
onto the second which requires minimum cost.
\begin{rmk}
	In our framework, the mass to be moved corresponds to that of vehicles. We therefore measure the distance between two traffic density distributions by computing the minimal cost to rearrange vehicles from the first configuration onto the second one.
	We also stress that, since we are considering macroscopic models, vehicles are indistinguishable so we are not able to pinpoint the exact same vehicle in the two scenarios, cf.\ \cite{cristiani2019M2AS}. 
\end{rmk}
%


In the particular case of 1D optimal transport  various characterizations give alternative, more manageable, definitions of WD. 
For example, 
if $(X,\mathfrak{D})$ is
$\mathbb R$ with the euclidean distance $|\cdot|$,
then, in the case of two Dirac delta functions $\delta_x$, $\delta_y$
we simply have $W_p(\delta_x,\delta_y)=\mathfrak D(x,y):=|y-x|$.
That coincides with the
desired distance quantification in the case of only two vehicles on a road, located in $x$ and $y$ respectively. 
Instead, 
in the case of
two absolutely continuous measures with $d\muuno=\rhouno dx$ and $d\mudue=\rhodue dx$, 
by~\cite[Remark 2.19]{villani2003book} we have
\begin{equation}\label{caratterizzazioneW1inR}
	W_1(\muuno,\mudue)=
	\int_{\R}|\Funo(x)-\Fdue(x)|dx=
	\int_0^1\big|(\Funo)^{-1}(t)-(\Fdue)^{-1}(t)\big|dt, 
\end{equation}
with
\begin{equation}\label{caratterizzazioneW1inRbis}
	\Funo(x):=\int_{-\infty}^{x}\rhouno(z)dz, \quad \Fdue(x):=\int_{-\infty}^{x}\rhodue(z)dz
\end{equation}
for $p=1$, whereas for $p\geq 1$,
\begin{equation}\label{caratterizzazioneWpinR}
	W_p(\muuno,\mudue)=\left(\int_{\R}|T^*(x)-x|^p\ \rhouno(x) dx\right)^{1/p},
\end{equation}
where $T^*:\R\to\R$ satisfies
\begin{equation}\label{caratterizzazioneWpinRbis}
	\int_{-\infty}^{x} \rhouno(z)dz=
	\int_{-\infty}^{T^*(x)} \rhodue(z)dz
	\quad \forall x\in\mathbb R.
\end{equation}

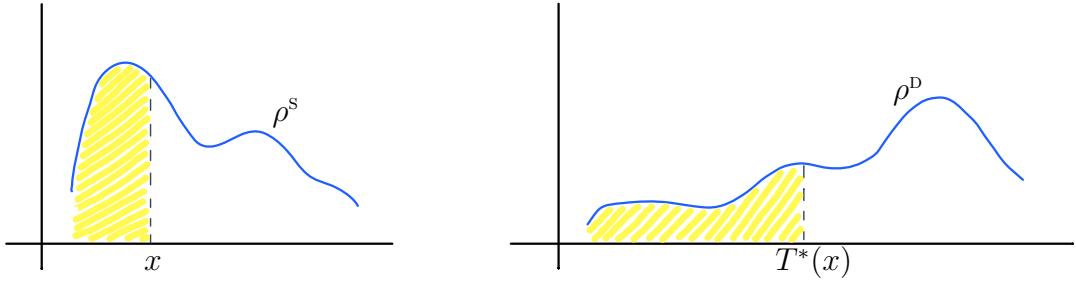
\begin{figure}[h!]
	\centering
	\begin{overpic}[width=0.8\textwidth]{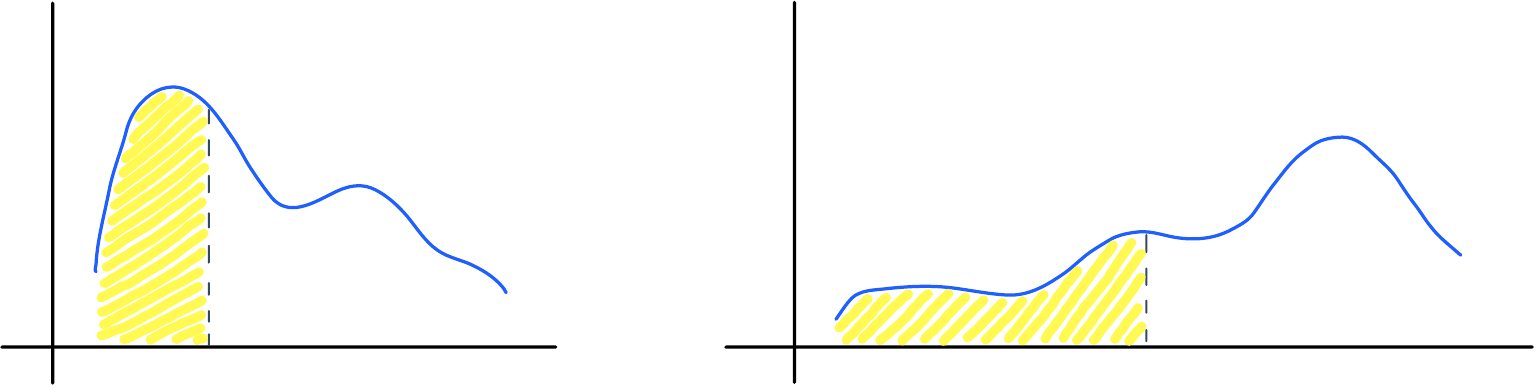}
		\put(13,0){$x$}
		\put(72,0){$T^*(x)$}
		\put(25,14){$\rhouno$}
		\put(83,16){$\rhodue$}
	\end{overpic}
	\caption{Property of monotone rearrangement.}
	\label{fig:monotonerearrangement}
\end{figure}

Equations \eqref{caratterizzazioneWpinR}-\eqref{caratterizzazioneWpinRbis} translate the fundamental property of \emph{monotone rearrangement} of the mass \cite[Remark 2.19]{villani2003book}, i.e.\ the optimal strategy consists in transferring the mass starting from the left, see Figure \ref{fig:monotonerearrangement}.

\subsection{Numerical setting}\label{sec:WD-numerica}
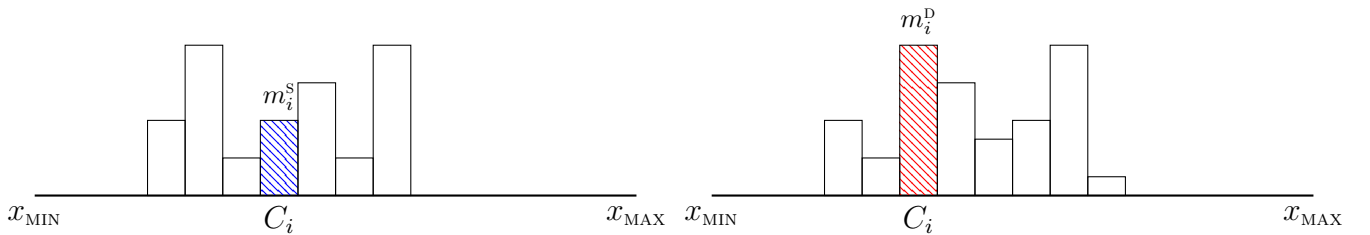
\begin{figure}[h!]
	\centering
	\begin{tikzpicture}[scale=1.0]
	\draw[thick] (-8.5,0) -- (-0.5,0);
	\node at (-8.5,0) [below] {$\xmin$};
	\node at (-0.5,0) [below] {$\xmax$};
	\draw[thick] (0.5,0) -- (8.5,0);
	\node at (0.5,0) [below] {$\xmin$};
	\node at (8.5,0) [below] {$\xmax$};
	\draw (-7,0) rectangle (-6.5,1);
	\draw (-6.5,0) rectangle (-6,2);
	\draw (-6,0) rectangle (-5.5,0.5);
	\draw[pattern=north west lines, pattern color=blue] (-5.5,0) rectangle (-5,1); \node at (-5.25,0)[below]{$C_i$}; \node at (-5.25,1)[above]{\footnotesize $m^\supply_i$};
	\draw (-5,0) rectangle (-4.5,1.5);
	\draw (-4.5,0) rectangle (-4,0.5);
	\draw (-4,0) rectangle (-3.5,2);
	\draw (2,0) rectangle (2.5,1);
	\draw (2.5,0) rectangle (3,0.5);
	\draw[pattern=north west lines, pattern color=red] (3,0) rectangle (3.5,2); \node at (3.25,0) [below] {$C_i$}; \node at (3.25,2)[above]{\footnotesize$m^\demand_i$};
	\draw (3.5,0) rectangle (4,1.5);
	\draw (4,0) rectangle (4.5,0.75);
	\draw (4.5,0) rectangle (5,1);
	\draw (5,0) rectangle (5.5,2);
	\draw (5.5,0) rectangle (6,0.25);
\end{tikzpicture}
	\caption{Numerical setting for a generic one-dimensional problem. Two copies of the domain $X$ with supply mass only (left) and demand mass only (right).}
	\label{fig:numericalsetting}
\end{figure}

Let us introduce here the numerical setting that will accompany us throughout the paper, see Figure \ref{fig:numericalsetting}.
For simplicity, we will only consider the one-dimensional interval $X=[\xmin,\xmax]$, divided in $N$ equispaced cells $C_1,\ldots,C_N$
($N$ is our first parameter of discretization, and we expect to vanish the numerical error for $N\to +\infty$).
We also denote by $x_i$ the barycenter of the cell $C_i$, $i=1,\ldots,N$ and by $\Dx:=x_{i+1}-x_i$ the step of the space discretization.
To measure the distance between two cells $C_i$ and $C_j$, we will always use the Euclidean norm $\mathfrak D(x_j,x_k):=|x_k-x_j|$.

Moving to the masses, for any $i=1,\ldots,N$ we denote by 
\begin{equation}\label{from-rho-to-m}
	m^\supply_i:=\rho^\supply(x_i)\Dx, \qquad 
	m^\demand_i:=\rho^\demand(x_i)\Dx, 
\end{equation} 
the supply and the demand mass located in $C_i$, respectively.
Roughly speaking, this corresponds to concentrate the masses on specific points, thus getting purely atomic density distributions 
$\hat\rho^\supply:=\sum_i m_i^\supply\delta_{x_i}$ 
and 
$\hat\rho^\demand:=\sum_i m_i^\demand\delta_{x_i}$
where $\delta$ denotes the Dirac delta.
Note that, in our context, the vectors $m^\supply:=\{m_i^\supply\}_i$ and $m^\demand:=\{m_i^\demand\}_i$ are given.
Moreover, to avoid any confusion between theoretical and numerical boundary conditions we will always assume that 
\begin{equation}\label{CBnulle}
	m_1^\supply=m_N^\supply=m_1^\demand=m_N^\demand=0,
\end{equation}
in such a way that boundaries remain out of play.

For some reasons which will be clear later on, it is also useful defining the \emph{remaining supply mass at} $C_i$, denoted by $\mathring m^\supply_i$, and the \emph{remaining demand mass} at $C_i$, denoted by $\mathring m^\demand_i$, $i=1,\ldots,N$.

\medskip

In order to get a numerical approximation of the WD in $\R^1$, it is convenient to resort to the definition \eqref{caratterizzazioneW1inR}. 
Indeed, any quadrature formula is suitable to get the desired value with arbitrary precision, provided $\rho^{\supply},\rho^{\demand}$ are sufficiently smooth and $N$ is sufficiently large.
For example, the naive composite mid-point formula gives
\begin{equation}\label{WD-numerica-rettangoli}
	W_1(\rho^\supply,\rho^\demand)\approx
	\tilde W_1(m^\supply,m^\demand):=\Dx
	\sum_{j=1}^{N}
	\left| 
	\sum_{i=1}^{j}m_i^\supply -
	\sum_{i=1}^{j}m_i^\demand
	\right|.
\end{equation}

The case of more general spaces, including networks embedded in $\R^2$, is more difficult. Typically the problem is solved by reformulating it in a fully discrete setting as a Linear Programming (LP) problem. 
The interested reader can find in the books by Santambrogio \cite[Section 6.4.1]{santambrogiobook} and Sinha \cite[Chapt.\ 19]{sinhabook} the complete procedure, recently used also in the context of traffic flow on networks \cite{briani2018CMS}.
Even if we will stick with the one-dimensional case, we recall now the crucial steps of this approach because it will be useful in the following. 

\medskip

Following the mass transport interpretation, let $c_{jk}$ be the cost of shipping a unit quantity of mass from the source $C_j$, $j=1,\ldots,N$, to the destination $C_k$, $k=1,\ldots,N$. 
In general, $c_{jk}$ is usually defined as the length of the shortest path joining $x_j$ and $x_k$. 
Let $\gamma_{jk}$ be the (unknown) quantity shipped from the source $j$ to the destination $k$. 
The problem is then formulated as
\begin{equation}\label{hitchcock}
	\begin{tabular}{ll}
		minimize & $\sum\limits_{j=1}^N \sum\limits_{k=1}^N c_{jk} \gamma_{jk}$ \\ [4mm]
		subject to & $\sum\limits_k \gamma_{jk}=m^\supply_j,\quad \forall j$ \\  [3mm]
		& $\sum\limits_j 
		\gamma_{jk}=m^\demand_k,\quad \forall k$ \\  [2mm]
		& $\gamma_{jk}\geq 0$.
	\end{tabular}
\end{equation}
The solution gives an approximation $\tilde W_1$ of $W_1$.
Note that the solution satisfies $\gamma_{jk}\leq\min\{m^\supply_j,m^\demand_k\}$ since one cannot move more than $m^\supply_j$ from any source $j$ and it is useless to bring more than $m^\demand_k$ to any destination $k$.
From \eqref{hitchcock}, it is easy to recover a standard LP problem with equality constraints
\begin{equation}\label{PL}
	\begin{tabular}{ll}
		minimize & $\mathbf{c}^\intercal \mathbf{x}$ \\ [3mm]
		subject to & $\mathbf{A}\mathbf{x}=\mathbf{b}$ \\  [2mm]
		& $\mathbf{x}\geq 0$,
	\end{tabular}
\end{equation}
simply defining
\begin{equation}\label{def:x,c,b}
	\begin{array}{l}
		\mathbf{x}:=(\gamma_{11},\gamma_{12},\ldots,\gamma_{1N},\gamma_{21},\gamma_{22},\ldots,\gamma_{2N},\ldots,\gamma_{N 1},\ldots,\gamma_{N N})^\intercal \\
		\mathbf{c}:=(c_{11},c_{12},\ldots,c_{1N},c_{21},c_{22},\ldots,c_{2N},\ldots,c_{N 1},\ldots,c_{N N})^\intercal \\
		\mathbf{b}:=(m^\supply_1,\ldots,m^\supply_{N},m^\demand_1,\ldots,m^\demand_{N})^\intercal 
	\end{array}
\end{equation}
and $\mathbf{A}$ as the $2N\times N^2$ block matrix
\begin{equation}\label{def:matrixA}
	\mathbf{A}:=
	\left[
	\begin{array}{ccccc}
		\unoN & \zeroN & \zeroN & \cdots & \zeroN \\
		\zeroN & \unoN & \zeroN & \cdots & \zeroN \\
		\zeroN & \zeroN & \unoN & \cdots & \zeroN \\
		\vdots & \vdots & \vdots & \ddots & \vdots \\
		\zeroN & \zeroN & \zeroN & \cdots & \unoN \\
		I_N & I_N & I_N & I_N & I_N \\
	\end{array}
	\right],
\end{equation}
where $I_N$ is the $N\times N$ identity matrix, 
$\unoN:=(\underbrace{1\ 1\ \cdots \ 1}_{N \text{ times}})$, 
and $\zeroN:=(\underbrace{0\ 0\ \cdots \ 0}_{N \text{ times}})$.

\section{Basic theory for Generalized Wasserstein distances}\label{sec:GWD}
In this section we recall some definitions of metrics
on the space of positive finite Borel
measures on the locally compact
complete and separable metric
space $(X,\mathfrak D)$, 
inducing the weak-narrow convergence,
all inspired to the distance defined by \eqref{def:W}
between probability measures.

\subsection{Figalli \& Gigli approach \cite{figalli2010JMPA}}\label{sec:GWD-figalli-teoria}

In \cite{figalli2010JMPA}, an optimal transport type distance is introduced,
with respect to which the heat equation under Dirichlet conditions at the boundary
of an open bounded domain $\Omega$ can be described as the gradient
flow of an entropy, similarly to the case of the heat equation in free space. 
In order to encode the correct boundary conditions in the minimizing movement scheme, the 
GWD between  two finite non-negative (Borel)
measures $\mu^{\supply},\mu^{\demand}$ on $\Omega$ is defined by 
\begin{equation}
	\label{FGW}
	W^{\fg}_p(\mu^{\supply},\mu^{\demand}) = 
	\left\{\inf_{\gamma} \int_{\,\overline\Omega\times \overline\Omega} \|x-y\|^p d\gamma(x,y)\right\}^\frac1p ,
\end{equation}
where infimum is over
\begin{equation}
	\label{ADM-FG}
	\Big\{
	\gamma\in \mathcal{M}^+(\overline\Omega\times\overline\Omega)\
	\Big\vert\
	\big(\pi^\supply_{\,\,\sharp} \gamma\big)\mres\Omega=\mu^{\supply} \,,\
	\big(\pi^\demand_{\,\,\sharp} \gamma\big)\mres\Omega=\mu^{\demand} 
	\Big\}\,.
\end{equation}
In the latter, we used the notation $(\bar{\mu}\mres\Omega)(A)=\bar{\mu}(A\cap\Omega) $
for the restriction to $\Omega$ of a measure $\bar{\mu}$ on $\overline\Omega$, and 
for all given $\gamma$
we denoted by $\pi^\supply_{\,\,\sharp} \gamma$ (resp., $\pi^\demand_{\,\,\sharp} \gamma$)
the push-forward of $\gamma$ under the projection defined on $\overline\Omega\times\overline\Omega$ by
$\pi^\supply(x,y)=x$ (resp., $\pi^\demand(x,y)=y$). Equivalently, a non-negative
finite Radon measure $\gamma$ on $\overline\Omega\times\overline\Omega$ is
an admissible competitor for the minimization in \eqref{FGW} if and only if
\begin{equation}
	\gamma(A \times \overline\Omega) = \mu^{\supply}(A)\,,\ \gamma(\overline\Omega\times B)=\mu^{\demand}(B)\,,\qquad \text{for all $A,B\in\mathcal{B}(\Omega)$.}
\end{equation}

Let us note that \eqref{ADM-FG} coincides with the set 
\[
\bigcup_{({\overline\mu}^{\supply},{\overline\mu}^{\demand})}\Big\{\gamma\in \mathcal{M}^+(\overline\Omega\times\overline\Omega)\ \Big\vert\
\pi^\supply_{\,\sharp}\gamma={\bar\mu}^{\supply}\,,\ 
\pi^\demand_{\,\sharp}\gamma={\bar\mu}^{\demand}
\Big\},
\]
where the union is over all pairs
$({\bar\mu}^{\supply},{\bar\mu}^{\demand})\in \mathcal{M}^+(\overline\Omega)\times
\mathcal{M}^+(\overline\Omega)$ such that
\[
{\bar\mu}^{\supply}\mres\Omega=\mu^{\supply}\,,\
{\bar\mu}^{\demand}\mres\Omega=\mu^{\demand}\,,\
{\bar\mu}^{\supply}(\overline\Omega)={\bar\mu}^{\demand}(\overline\Omega)\,.
\]
Hence, an equivalent variational definition of the unbalanced optimal transport distance
\eqref{FGW}
is
\begin{equation}
	\label{FGW2}
	\min\Big\{
	W_p({\bar\mu}^{\supply},{\bar\mu}^{\demand})\
	\Big\vert\
	{\bar\mu}^{\supply}\mres\Omega=\mu^{\supply}\,,\
	{\bar\mu}^{\demand}\mres\Omega=\mu^{\demand}\,,\
	{\bar\mu}^{\supply}(\partial\Omega)-{\bar\mu}^{\demand}(\partial\Omega)=
	\mu^{\supply}(\Omega)-\mu^{\demand}(\Omega)
	\Big\},
\end{equation}
where the classical WD is defined as in \eqref{def:W} with $\mathfrak D(x,y)=\|x-y\|$ for probability
measures, and extended by homogeneity to non-negative Radon measures with equal masses. With this definition, if $\bar\mu^\supply,\bar\mu^\demand\in\mathcal{M}^+(\overline\Omega)$ are
fixed, then $W_p(t \bar\mu^\supply ,t\bar\mu^\demand)$ is arbitrarily large
if $t>0$ is large enough. Hence the minimum
in \eqref{FGW2} is indeed achieved.
Notice that mass unbalance $\mu^{\supply}(\Omega)-\mu^{\demand}(\Omega)$
is fixed if so are the marginals $\mu^{\supply}$ and $\mu^{\demand}$. 

The authors of~\cite{figalli2010JMPA} deal explicitly with the case of a quadratic unitary cost, i.e., $c(x,y)=\|x-y\|^2$ in \eqref{FGW}, but their main results also hold for
a more general cost, e.g.,
$c(x,y)=\|x-y\|^p$, with $p\ge1$. As intermediate steps towards their main results,
they prove lower semicontinuity with respect to the weak-$\ast$ convergence on $\Omega$ and provide a gluing lemma implying the triangle inequality, so that \eqref{FGW} defines
a distance on $\mathcal{M}^+(\Omega)$, and the infimum is achieved by an
appropriate optimal plan $\gamma$.

The idea of minimizing \eqref{FGW} among all 
$\gamma$ as in \eqref{ADM-FG}
is that $\partial\Omega$ models an infinite
mass reservoir and the transportation policy described by
convenient competitors $\gamma$  might include the option to cut out part of the material from the actual transference plan, by transferring it to $\partial\Omega$, as well as the possibility to integrate the transportation plan by transferring material from $\partial\Omega$ to the final destination, provided that the transportation cost to and from $\partial\Omega$ is accounted for in the overall integral cost together with that
of moving the mass that is authentically
transported.
A dual formulation of the problem leads to a generalization of Kantorovich potentials and 
Kantorovich norms, which had been already investigated especially in the case of a linear cost ($p=1$),
see~\cite{Hanin1992PAMS, KantorovichRubinstein1958VestnLeningrUniv}.

\begin{rmk}\label{FGnotext}
	It is clear from \eqref{FGW2} that the distance
	between positive finite Radon measures introduced
	in \cite{figalli2010JMPA} is not an extension of the
	classical Wasserstein distance: if the two prescribed marginals are concentrated close to the boundary with supports far apart from each other, then
	while solving \eqref{FGW2} there will be a convenience
	for $\bar\mu^\supply$ and $\bar\mu^\demand$ in charging
	appropriately two corresponding disjoint boundary regions, after which the Wasserstein to be minimized will be achieved by transport plans that only move mass from or to the boundary. In fact, the goal of \cite{figalli2010JMPA} was to generalize Otto's calculus to PDEs under Dirichlet conditions, rather than extending the optimal transport metrics to 
	finite positive Radon measures.
\end{rmk}

\subsection{Piccoli \& Rossi approach \cite{piccoli2014ARMA}}\label{sec:GWD-piccoli-teoria}
In \cite{piccoli2014ARMA}, the set
of positive finite Borel measures on $\Omega$
is endowed with a metrics, defined by

\begin{equation}\label{def:GWD_PR}
	W_{p}^\pr(\mis{\supply},\mis{\demand};a,b):=
	\inf_\gamma 
	\Bigg\{
	a\big\vert\mis{\supply}(\Omega)-\pi^\supply_{\,\sharp}\gamma (\Omega)\big\vert
	+
	a\big\vert\mis{\demand}(\Omega)-
	\pi^\demand_{\, \sharp}\gamma(\Omega)\big\vert
	+
	b\int_{\Omega\times\Omega} \|x-y\|^p\,d\gamma 
	\Bigg\}^\frac1p,
\end{equation}
where $a,b>0$ are fixed parameters and the infimum
is among all finite positive Borel measures
$\gamma$ on $\Omega\times\Omega$.

Note that \eqref{def:GWD_PR} can be equivalently written in the form
\begin{equation}\label{def:GWD_PR2}
	W_{p}^\pr(\mis{\supply},\mis{\demand};a,b)=
	\inf 
	\Big\{
	a\big(
	\big\vert\mis{\supply}(\Omega)-\miss{\supply}(\Omega)\big\vert
	+
	\big\vert\mis{\demand}(\Omega)-\miss{\demand}(\Omega)\big\vert
	\big)
	+
	bW_p\big(\miss{\supply},\miss{\demand}\big)
	\ \Big\vert\ 
	\miss{\supply}(\Omega)=\miss{\demand}(\Omega)
	\Big\},
\end{equation}
where $\mis{\supply}$ and $\mis{\demand} $ form an
arbitrary pair of positive finite Borel measure on $\Omega$ with equal mass. In the right hand-side of 
\eqref{def:GWD_PR2},
the term $W_p(\miss{\supply},\miss{\demand})$ is understood as the value taken by the extension
by homogeneity to pairs of finite positive measures of equal mass of the Wasserstein distance
between probability measures.\footnote{
	Namely, we are setting
	$W_p(\miss{\supply},\miss{\demand})=
	m W_p(\miss{\supply}/m,\miss{\demand}/m)$, where
	$m=\miss{\demand}(\Omega)=
	\miss{\supply}(\Omega)$ and $\miss{\supply},\miss{\demand}$
	are thought of as measures on $\overline\Omega$ that
	vanish on all measurable subsets of the boundary
	$\partial\Omega=\overline{\Omega}\setminus\Omega$.
}

The infimum in \eqref{def:GWD_PR2}, and hence in \eqref{def:GWD_PR} too,
is indeed achieved and it defines a
distance~\cite[Prop.\ 1]{piccoli2014ARMA}.

The mass transfer may include either
physical transportation or creation and destruction. Adding mass to and removing it from $\mis{\supply}$, as well as from $\mis{\demand}$, are both allowed actions, with the implied costs being weighted by the parameter $a$. After these actions have been taken,
the supplies eligible for transport,
as well as the actual demands, are established and
the corresponding optimal transfer policy
contributes to the overall cost, via a classical WD weighted by 
$b$.

The constraints
$\miss{\supply}(\Omega)\le \mu^\supply(\Omega)$
and $\miss{\demand}(\Omega)\le \mu^\demand(\Omega)$ may be added in the minimization
without changing the minimum value, see~\cite[Prop.\ 2]{piccoli2016ARMA}.

\begin{rmk}
	\label{PRnotext}
	Heuristically, all
	combinations of the above actions are admissible for the minimization:
	an arbitrary 
	fraction of mass may be removed from the marginals, along with transporting  the remaining mass. Accordingly, \eqref{def:GWD_PR}
	does not extend the
	Wasserstein distance because in the case
	of two marginals with equal masses the allowed actions
	include the possibility of removing mass from both.
	On the other hand, in the case $p=1$ with $a=b=1$
	it still offers a natural generalization of the 
	$1$-WD, being equivalent to
	the flat metric.
\end{rmk}

The Piccoli \& Rossi distance
\eqref{def:GWD_PR} differs from the Figalli \& Gigli one \eqref{FGW}. 
Recalling \eqref{FGW2}, it is clear that,
in the latter, the mass added to the marginals is located at the boundary and that no other pay-off is required except for the transportation cost between all locations, including those at the boundary. 
By contrast, the boundary geometry plays no role in \eqref{def:GWD_PR}.

On the other hand, at least in the case of 
a quadratic cost ($p=2$), it is known (see \cite[Prop.\ 1.1]{ChiPeySchVia2018FoundCompMath}) that
\eqref{def:GWD_PR} is equivalent to the partial
optimal transport problem~\cite{CafMcC2010ANNALS,Figalli2010ARMA}, i.e.,
\[
\inf_\gamma \int_{\Omega\times\Omega}
\Big( \|x-y\|^2-\lambda \Big)\,d\gamma \,,
\]
where the infimum is among all non-negative
Radon measures on $\Omega\times\Omega$,
whose left and right marginals are dominated by
$\mis{\supply} $ and $\mis{\demand}$,
that are subject to have
a fixed overall mass $m$, and the Lagrange multiplier
for this constraint satisfies $\lambda=2 a^2/b^2$.

Furthermore, the distance \eqref{def:GWD_PR} is inspired
by the dynamical formulation of the Monge--Kantorovich problem: by generalizing the Benamou and Brenier's
approach~\cite{BenBre2000NumerMath}, the
supply and the demand distributions are thought of
as the initial and the final prescribed configuration
for a fluid motion and the 
minimal value of \eqref{def:GWD_PR}
is obtained by minimizing the total kinetic energy
necessary to connect the two configurations so as to
balance mass conservation with
creation and destruction of mass, thanks to the presence of an additional source term in a differential constraint for the macroscopic fluid variables.

\subsection{Gaussian Hellinger--Kantorovich approach~\cite{LieMieSav2018InventMath}}\label{sec:GWD-savare-teoria}
Another distance between positive finite Borel measures on $\Omega$ is defined by
\begin{subequations}\label{def:LET-HK}
	\begin{equation}
		\label{def:HKdist}
		W^\ghk_{2}(\mu^\supply,\mu^\demand;a,b)
		:=
		\inf
		\Bigg\{
		a\,\mathcal{U}\big( 
		\pi^\supply_{\,\sharp}\gamma \mathbin\vert \mu^\supply
		\big)
		+
		a\,\mathcal{U}\big( 
		\pi^\demand_{\,\sharp}\gamma \mathbin\vert \mu^\demand
		\big)
		+
		b\int_{\Omega\times\Omega}
		\|x-y\|^2\,d\gamma
		\Bigg\}^\frac{1}{2},
	\end{equation}
	with $a,b>0$ being fixed parameters, and the infimum is among all
	positive finite measures $\gamma$ on $\Omega\times\Omega$. 
	Here,
	\begin{equation}\label{LET}
		\mathcal{U}\big(
		\tilde{\mu}     \mathbin\vert \mu
		\big)=
		\begin{cases}
			\displaystyle
			\int U(\sigma)\,d\mu \,, & 
			\text{if $d\tilde{\mu}=\sigma d \mu$}\,,\\
			+\infty \,, & 
			\text{otherwise,}
		\end{cases}\,,\qquad\text{where }
		U(\sigma) = \sigma \log(\sigma) -\sigma+1\,.
	\end{equation}
\end{subequations}

Such a distance 
is called
in~\cite{LieMieSav2018InventMath}
the {\em Gaussian Hellinger--Kantorovich distance} in view of the occurrence of a Gaussian function
in its homogeneous formulation, see Section \ref{sec:SS-teoria} below.
We note that \eqref{def:LET-HK} can be equivalently written in the form
\begin{equation}
	\label{def:HKdist2}
	W^\ghk_{2}(\mu^\supply,\mu^\demand;a,b)
	=
	\inf
	\Bigg\{
	a\,\mathcal{U}\big( 
	\hat{\mu}^\supply \mathbin\vert \mu^\supply 
	\big)
	+
	a\,\mathcal{U}\big( 
	\hat{\mu}^\demand \mathbin\vert \mu^\demand
	\big)
	+
	b W_2^2(\hat{\mu}^\supply,\hat{\mu}^\demand)
	\Bigg\}^\frac{1}{2},
\end{equation}
where the infimum is among all
positive finite measures $\hat{\mu}^\supply$
and $\hat{\mu}^\demand$ of equal mass on $\Omega$.
In order to define the
last term in the right hand-side of \eqref{def:HKdist2},
we have implicitly understood 
both $\hat{\mu}^\supply$ and $\hat{\mu}^\demand$
as multiples of probability measures
on $\overline{\Omega} $ vanishing on $\partial\Omega$
and we have defined their WD accordingly,
with the obvious extension by homogeneity
to pairs of finite positive measures of equal mass on $\Omega$
of the WD between probability
measures.

When the given measures $\mu^\supply$ and
$\mu^\demand$ themselves happen to have equal mass, 
they form an admissible pair for the minimization: then, the right hand side of \eqref{def:HKdist2}
is the classical (quadratic) WD, whose
the value is a computed by solving a balanced optimal transport problem, because in view of \eqref{LET}
we have $U(1)=0$.
In the unbalanced case, the two additional
non-negative contributions
in the objective penalize
the deviation of the competitors from the prescribed
measures.

The variational principle \eqref{def:LET-HK} belongs to the 
wide class of {\em optimal entropy-transport problems} considered 
in~\cite{LieMieSav2018InventMath}, in which
the penalization
is proportional to the {\em Csisz\'ar--Morimoto
	divergence}, defined in terms of
proper convex functions $U:[0,+\infty)\to [0,+\infty]$,
called the {\em entropy functions},
possibly different from \eqref{LET}.
That large class covers the
Piccoli \& Rossi distance \eqref{def:GWD_PR2} for $p=1$, 
with entropy function $U(\sigma)=|\sigma-1|$,
and is comprised of
several models, alternative to~\cite{piccoli2014ARMA}
but equally inspired by dynamical formulations with source terms in the differential constraints, which have appeared in the literature since.

In particular, the model case 
of the {\em logarithmic optimal entropy-transport problem}
\eqref{def:LET-HK}
has been investigated under various perspectives:
its value coincides with that of the metrics studied 
in~\cite{ChiPeySchVia2018FoundCompMath, KonMorVor2016AdvDE}
(also called Wasserstein--Fisher--Rao distance, as it interpolates between the quadratic Wasserstein
and the Fisher--Rao distance), for which a static definition had been provided in~\cite{ChiPeySchVia2018JFA}.

\begin{rmk}\label{GHKsiext}
	It is known, e.g. by~\cite[Part~II]{LieMieSav2018InventMath},
	that
	the logarithmic entropy-transport minimization
	\eqref{def:LET-HK}
	defines a distance 
	that metrizes the weak/narrow convergence of
	non-negative finite measure,
	extending the quadratic WD 
	(otherwise defined only for probability measures).
\end{rmk}

\subsection{Savar\'e \& Sodini approach \cite{savare2024JMPA}}\label{sec:SS-teoria}
In~\cite[Part~II]{LieMieSav2018InventMath} 
the logarithmic entropy-transport minimization
\eqref{def:LET-HK}
is proved to be equivalent to an optimal transport type problem in the 
{\em metric cone} $\mathfrak{C}[X]$ over $X$.
More precisely, one defines $\mathfrak{C}[X]$ as the quotient of 
$X\times[0,+\infty)$
obtained by collapsing the whole fiber $X\times \{0\}$
to a single element, endowed with a distance that is
defined, once an appropriate function $g$ is chosen, by
setting 
\[
d_{\mathfrak{C},g}([x^\supply,r^\supply],[x^\demand,r^\demand])=
\sqrt{H_g(x^\supply,r^\supply;x^\demand,r^\demand)}\,,
\]
where
\begin{equation}
	\label{def:H-SS}
	H_g(x,r;y,s)=r+s -2\sqrt{rs}\ g(\|x-y\|)
	\quad\text{and}
	\quad g(t)= e^{-\frac{t^2}{2}}\,.
\end{equation}
Then, one considers
the minimization of
\[
\int_{\mathfrak{C}[X]\times \mathfrak{C}[X]} 
d^2_{\mathfrak{C},g}([x^\supply,r^\supply],
[x^\demand,r^\demand]) \,d\alpha
\]
among positive measures $\alpha$ on
$\mathfrak{C}[X]\times\mathfrak{C}[X]$
with {\em homogeneous marginals} $\mu^\supply$ and $\mu^\demand$, i.e., such that
\[
\mu^\supply(A)=
\int_{\mathfrak{C}[A]\times\mathfrak{C}[X]} r^\supply\,d\alpha([x^\supply,r^\supply],[x^\demand,r^\demand])\,,\quad
\mu^\demand(B)=
\int_{\mathfrak{C}[X]\times\mathfrak{C}[B]} 
r^\demand\,d\alpha([x^\supply,r^\supply],
[x^\demand,r^\demand])
\]
for all Borel sets $A,B\subset X$.

The constrained minimum is achieved, and
it coincides with the 
Gaussian Hellinger--Kantorovich 
distance defined by \eqref{def:LET-HK} in the case $a=b=1$.\footnote{
	A match is  possible also for different choices of parameters $a,b$, provided that the definition of $H_g$ changes accordingly in \eqref{def:H-SS}. 
	Moreover, the choice $g(t)=\cos\left( t^2 \wedge \tfrac{\pi}{2}\right)$
	is also possible in \eqref{def:H-SS}, leading to the {\em Hellinger--Kantorovich distance}, which equals \eqref{def:HKdist} if \eqref{LET} is replaced by the entropy function
	$U(\sigma)=\sigma-1-\log\sigma $.
}

\section{Discrete and computational approaches for Generalized Wasserstein distances}\label{sec:GWDcomputational}
In this section, we first revise the GWDs previously introduced in the special case of atomic measures (sum of Dirac delta's). 
This allows us to naturally introduce a discretization in the space domain, which is the first step towards a numerical computation. 
Secondly, we propose original computational approaches to get algorithms ready to be implemented on a computer. 
In the final Section \ref{sec:comparison}, we present a comparative numerical test of all the proposed approaches. 
Specific tests for each GWD are instead postponed in Appendix \ref{appendice}.

\subsection{Figalli \& Gigli approach}\label{sec:GWD-figalli-computational}

\subsubsection{Discrete approach}\label{sec:GWD-figalli-discrete}
In the special case of two purely atomic marginals with finite supports in a one-dimensional domain, i.e.,
\[
\mu^\supply = \sum_{i=2}^{N-1} m_i^\supply\delta_{x_i}\,,\qquad
\mu^\demand = \sum_{i=2}^{N-1} m_i^\demand \delta_{x_i}\,,
\]
we get an interpretation of
\eqref{FGW2}
as a linear programming problem (indices 1 and $N$ are discarded because of \eqref{CBnulle}). 
Namely,
we consider the problem
\[
\min_{(m_1^\supply,m_N^\supply,m_1^\demand,m_N^\demand)\in\mathbb R^4}\left\{
W_p^p\Big( 
\sum_{i=1}^{N} m_i^\supply\delta_{x_i}
\mathbin,
\sum_{i=1}^{N} m_i^\demand\delta_{x_i}
\Big)\ \Big\vert\ 
m_1^\supply+m_N^\supply -(m_1^\demand+m_N^\demand)=
\Delta M
\right\}\,,
\]
where 
$$
\Delta M:=\sum_{i=2}^{N-1}m_i^\supply-\sum_{i=2}^{N-1}m_i^\demand,
$$
and, after setting $c_{jk}=\|x_j-x_k\|^p$, 
for all
$j,k\in\{1,\ldots,N\}$, we write it in the equivalent form 
\[
\begin{split}
	\tilde W^\fg_p(m^\supply,m^\demand):= 
	\min\Bigg\{
	\sum_{j,k=1}^N c_{jk}\gamma_{jk}\ \Big\vert\ 
	\gamma\in [0, +\infty)^{N\times N}\,,\  
	& \sum_{k=1}^N\gamma_{jk}=m_j^\supply\,,\
	\text{for $j\in\{2,\ldots,N-1\}$}\,,\\
	& \sum_{j=1}^N\gamma_{jk}=m_k^\demand\,,\
	\text{for $k\in\{2,\ldots,N-1\}$}, \\ 
	&  \sum_{k=1}^N\big( \gamma_{1k}+\gamma_{Nk}\big)
	-\sum_{j=1}^N\big(\gamma_{j1}+\gamma_{jN}\big)
	=
	\Delta M
	\Bigg\}^{\frac1p}\,.
\end{split}
\]
That amounts to the minimization of a linear objective
in a convex polytope obtained by
intersecting $2N-3$ affine half-spaces in $\mathbb R^{N\times N}$ with $[0,+\infty)^{N\times N}$. 
Alternatively, the last constraint
may be converted into an additional contribution to the total cost involving suitable
Lagrange multipliers.

\subsubsection{Computational approach}\label{sec:GWD-figalli-numerica}

The minimization problem introduced in the previous section can be formulated as a LP problem in the standard form \eqref{PL} defining
\begin{eqnarray*}
	& &\mathbf{x}:=(\gamma_{11},\gamma_{12},\ldots,\gamma_{1N},\gamma_{21},\gamma_{22},\ldots,\gamma_{2N},\ldots,\gamma_{N 1},\ldots,\gamma_{N N})^\intercal \\
	& &\mathbf{c}:=(c_{11},c_{12},\ldots,c_{1N},c_{21},c_{22},\ldots,c_{2N},\ldots,c_{N 1},\ldots,c_{N N})^\intercal \\
	& &\mathbf{b}:=(m^\supply_2,\ldots,m^\supply_{N-1},m^\demand_2,\ldots,m^\demand_{N-1},\Delta M)^\intercal
\end{eqnarray*}
and $\mathbf{A}$ as the $(2N-3)\times N^2$ block matrix
\begin{equation}\label{def:matrixA-FG-franz}
	\mathbf{A}:=
	\left[
	\begin{array}{cccccc}
		\zeroN & \unoN & \zeroN & \cdots & \zeroN & \zeroN \\
		\zeroN & \zeroN & \unoN & \cdots & \zeroN & \zeroN\\
		\vdots & \vdots & \vdots & \ddots & \vdots & \vdots\\
		\zeroN & \zeroN & \zeroN & \cdots & \unoN & \zeroN\\
		J_N & J_N & J_N & \cdots & J_N & J_N\\
		U_N & Z_N & Z_N & \cdots & Z_N & U_N
	\end{array}
	\right],
\end{equation}
where $I_N$ is the $N\times N$ identity matrix, 
$\unoN:=(\underbrace{1\ 1\ \cdots \ 1}_{N \text{ times}})$, 
$\zeroN:=(\underbrace{0\ 0\ \cdots \ 0}_{N \text{ times}})$,
$Z_N:=(-1\ \underbrace{0\ 0\ \cdots \ 0}_{N-2 \text{ times}} -1)$,
$U_N:=(0\ \underbrace{1\ 1\ \cdots \ 1}_{N-2 \text{ times}}\ 0)$,
and $J_N$ is the $(N-2)\times N$ block matrix $[\mathbb{0}_{\!N-2}^\intercal\ I_{N-2}\ \mathbb{0}_{\!N-2}^\intercal]$.

\medskip 
The LP formulation described above can be slightly generalized to account for the costs of creating/destroying mass at boundaries. 
Still keeping the simplifying assumption \eqref{CBnulle}, the vector $\mathbf x$ of the unknowns is modified adding four new unknowns $M_1^\supply, M_N^\supply, M_1^\demand, M_N^\demand$, corresponding to the masses taken from the left/right boundaries, and the masses brought to the left/right boundaries, respectively. 
Note that these masses are potentially unlimited. 
Moreover, four input parameters $c_1^\supply$, $c_N^\supply$, $c_1^\demand$, $c_N^\demand$ related to the four new unknowns are added. 
For example, the quantity $c_1^\supply M_1^\supply$ will give the cost to create the mass $M_1^\supply$ at the left boundary (N.B., just creating it, not transporting it inside the domain), while the quantity $c_N^\demand M_N^\demand$ will give the cost to destroy the mass $M_N^\demand$ at the right boundary.

The PL formulation becomes:
\begin{equation}\label{hitchcock_esteso}
	\begin{split}
		\tilde W^{\fg^\prime}_p(m^\supply,m^\demand):= 
		\min\Bigg\{ &
		\sum_{j,k=1}^N c_{jk}\gamma_{jk}+c_1^\supply M_1^\supply+c_N^\supply M_N^\supply+c_1^\demand M_1^\demand+c_N^\demand M_N^\demand \ \Big\vert\ 
		\gamma\in [0, +\infty)^{N\times N}\,,\\  
		& \sum_{k=1}^N\gamma_{jk}=m_j^\supply\,,\
		\text{for $j\in\{2,\ldots,N-1\}$}\,,\quad
		\sum_{j=1}^N\gamma_{jk}=m_k^\demand\,,\
		\text{for $k\in\{2,\ldots,N-1\}$}, \\ 
		&  \sum_{k=1}^N x_{jk}=M^\supply_j,\quad\text{for } j\in\{1,N\}, 
		\quad \sum_{j=1}^N x_{jk}=M^\demand_k,\quad\text{for } k\in\{1,N\}, \\
		& M_1^\supply, M_N^\supply, M_1^\demand, M_N^\demand\geq 0,
		\quad c_{11}=c_{1N}=c_{N1}=c_{NN}=+\infty
		\Bigg\}^{\frac1p}\,,
	\end{split}
\end{equation}


which corresponds to \eqref{PL} with
\begin{eqnarray*}
	& &\mathbf{x}:=(x_{11},x_{12},\ldots,x_{1N},x_{21},x_{22},\ldots,x_{2N},\ldots,x_{N 1},\ldots,x_{N N},
	M^\supply_1,M^\supply_N,M^\demand_1,M^\demand_N)^\intercal \\
	& &\mathbf{c}:=(c_{11},c_{12},\ldots,c_{1N},c_{21},c_{22},\ldots,c_{2N},\ldots,c_{N 1},\ldots,c_{N N},c_1^\supply,c_N^\supply,c_1^\demand,c_N^\demand)^\intercal \\
	& &\mathbf{b}:=(0,m^\supply_2,\ldots,m^\supply_{N-1},
	0,0,
	m^\demand_2,\ldots,m^\demand_{N-1},0)^\intercal
\end{eqnarray*}
and $\mathbf{A}$ is a $2N\times (N^2+4)$ block matrix, obtained from matrix \eqref{def:matrixA} by extending it with four columns. 
The extra columns have a single non-zero element each, namely 
\begin{equation}
	\mathbf{A}_{1,N^2+1}=-1,\quad 
	\mathbf{A}_{N,N^2+2}=-1, \quad 
	\mathbf{A}_{N+1,N^2+3}=-1, \quad 
	\mathbf{A}_{2N,N^2+4}=-1,
\end{equation}
due to the four extra unknowns and constraints of the problem \eqref{hitchcock_esteso}.
The last constraints of \eqref{hitchcock_esteso} are needed to assure that it never convenient creating mass at boundaries just to keep it there or to transport it on the other boundary. 


In Appendix \ref{sec:GWD-figalli-numericaltest}, we present two academic examples which serve both to illustrate the features of the Figalli \& Gigli approach and to facilitate a deeper understanding of the proposed methodology.
%
%
%
%
%
%
%
\subsection{Piccoli \& Rossi approach}\label{sec:GWD-piccoli-computational}
\subsubsection{Discrete approach}\label{sec:GWD-piccoli-discrete}
In the case of two purely atomic marginals
with finite supports
\[
\mis{\supply} = \sum_{i=1}^N m_i^\supply \delta_{x_i}\,,\qquad
\mis{\demand} = \sum_{i=1}^N m_i^\demand \delta_{x_i}\,,\qquad
\]
the distance \eqref{def:GWD_PR} is computed by
minimizing 
\[
a \sum_{i=1}^N
\big(m_i^\supply-\mathring m_i^\supply\big)
+ 
a  \sum_{i=1}^N 
\big(m_i^\demand-\mathring m_i^\demand\big)
+
b \sum_{j,k=1}^N c_{jk}\gamma_{jk}\,,
\]
where $c_{jk}=\|x_j-x_k\|^p$, 
among all $\gamma\in [0,+\infty)^{N\times N}$ 
such
that
\[
\begin{split}
	\sum_{k=1}^N
	\gamma_{jk}=\mathring m_j^\supply\,, 
	& \quad \text{for $j\in\{1,\ldots,N\}$},\\
	\sum_{j=1}^N\gamma_{jk}=\mathring m_k^\demand\,,
	& \quad \text{for $k\in\{1,\ldots,N\}$},
\end{split}
\]
and among all $\mathring m_1^\supply$,
\ldots, $\mathring m_N^\supply$,
$\mathring m_1^\demand$,
\ldots, $\mathring m_N^\demand$
with $\mathring m_i^\supply\le m_i^\supply$ and $\mathring m_i^\demand\le m_i^\demand$, $i\in\{1,\ldots,N\}$.
By taking into account such conditions, the objective
is to minimize
\[
a \sum_{i=1}^N\big(m_i^\supply+m_i^\demand\big) + 
\sum_{j,k=1}^N \big( b\, c_{jk} -2 a\big) \gamma_{jk}\,,
\]
among all $\gamma\in[0,+\infty)^{N\times N}$ such that
\[
\begin{split}
	&\sum_{k=1}^N\gamma_{jk}\le m_j^\supply\,,
	\qquad 
	\text{for all $j=1,\ldots,N$},\\
	&\sum_{j=1}^N\gamma_{jk}\le m_k^\demand\,,  
	\qquad 
	\text{for all $k=1,\ldots,N$}.
\end{split}
\]

\subsubsection{Computational approach}\label{sec:GWD-piccoli-numerica}
This time, the problem can be solved through a linear programming problem with inequality constraints. First, one solves the problem 
\begin{equation}\label{PLperPR}
	\begin{tabular}{ll}
		minimize & $\mathbf{\hat c}^\intercal \mathbf{x}$ \\ [3mm]
		subject to & $\mathbf{A}\mathbf{x}\leq\mathbf{b}$ \\  [2mm]
		& $\mathbf{x}\geq 0$,
	\end{tabular}
\end{equation}
where 
$
\mathbf{\hat c}:=b \mathbf c-2a \unoNq,
$
$\mathbf x$, $\mathbf c$, $\mathbf b$ are defined as in \eqref{def:x,c,b}, and $\mathbf A$ is defined as in \eqref{def:matrixA}. 
Then, once the solution $\mathbf x^*=\mathbf x^*(a,b)$ of \eqref{PLperPR} is found, the final GWD is given by 
$$\tilde W^\pr_1(m^\supply,m^\demand;a,b):=\mathbf x^*(a,b) + a \sum_{i=1}^N\big(m_i^\supply+m_i^\demand\big).
$$

In Appendix \ref{sec:GWD-piccoli-numericaltest}, we present two academic examples which serve both to illustrate the features of the Piccoli \& Rossi approach and to facilitate a deeper understanding of the proposed methodology.

%
%
%
%
%
%

\subsection{Gaussian Hellinger--Kantorovich approach}\label{sec:GWD-savare-computational}

\subsubsection{Discrete approach}\label{sec:GWD-savare-discrete}
When computing the distance \eqref{def:LET-HK}
of two purely atomic measures with finite supports
\[
\mu^\supply=\sum_{i=1}^Nm_i^\supply
\delta_{x_i}\,,
\qquad
\mu^\demand=\sum_{i=1}^Nm_i^\demand
\delta_{x_i}
\]
we may resort to a dual formulation. Namely,
by adapting the proof of \cite[Theorem 6.3]{LieMieSav2018InventMath} to the case
in which $a$, $b$ are not necessarily both equal to $1$, one sees that the distance equals
the maximum of
\[
a\sum_{i=1}^N m_i^{\supply} \delta_{x_i}\left(1-e^{-b\varphi_i/a}\right)
+
a\sum_{i=1}^N m_i^{\demand} \delta_{x_i}\left(1-e^{-b\psi_i/a}\right)
\]
among all $(\varphi_1,\ldots,\varphi_N,\psi_1,\ldots,\psi_N)\in\mathbb R^N\times\mathbb R^N $ such that
\begin{equation}\label{defcperGHK}
	\varphi_j+\psi_k\le c_{jk}:=|x_j-x_k|^2\,,\qquad\text{for all $j,k\in\{1,\ldots,N\}$.}
\end{equation}

\subsubsection{Computational approach}\label{sec:GHK-numerica}
In this case the optimization problem can be written as a nonlinear maximization problem with linear constraints. 
To this end, we define the function $f:\R^{2N}\to\R$, 
$$
f(\mathbf z):=
a\sum_{i=1}^N m_i^\supply \left(1-e^{-b z_i/a}\right) +
a\sum_{i=N+1}^{2N} m_{i-N}^\demand 
\left(1-e^{-b z_i/a}\right), \qquad \mathbf z\in\R^{2N}.
$$
Then, we aim at finding
$$\tilde W^\ghk_2(m^\supply,m^\demand):=
\max_{\mathbf z} f(\mathbf z)$$
under constraints
$$
\mathbf A^{\!\intercal} \mathbf z\leq \mathbf c,
$$
where $\mathbf c$ is defined as in \eqref{def:x,c,b} (with $c_{jk}$ as in \eqref{defcperGHK}) and $\mathbf A$ is the $N^2 \times 2N$ matrix defined as
$$
\mathbf A:=
\left[
\begin{array}{llll}
	\me_1 \ \me_1 \cdots \me_1 &
	\me_2 \ \me_2 \cdots \me_2 &
	\cdots &
	\me_N \ \me_N \cdots \me_N 
	\\
	\me_1 \ \me_2 \cdots \me_N &
	\me_1 \ \me_2 \cdots \me_N &
	\cdots &
	\me_1 \ \ \me_2 \ \cdots \me_N 
\end{array}
\right],
$$
with $\me_i=(0\ 0\ \cdots \underbrace{1}_{i^{th}} \cdots\ 0 \ 0)^\intercal\in\R^{N}$ being, as usual, the $i$th element of the canonical base of $\R^N$.

%
%
Numerical tests for the GWD discussed here are in common with those of the following GWD, see Appendix \ref{sec:GWD-SS-numericaltest}.    

%
%
%

\subsection{Savar\'e \& Sodini approach}\label{sec:SS-computational}

\subsubsection{Discrete approach}
For the purpose of empiric approximation, it may be convenient to resort to the convex
relaxation approach introduced in~\cite{savare2024JMPA}.
In particular, in view of~\cite[Theorem 3.5]{savare2024JMPA},
in the case of two atomic homogeneous marginals
\[
\mu^\supply=\sum_{i=1}^N m_i^\supply
\delta_{x_i}\,,
\qquad
\mu^\demand=\sum_{i=1}^N m_i^\demand
\delta_{x_i}
\]
the Gaussian Hellinger--Kantorovich distance can be obtained by minimizing
\begin{equation}\label{def:S}
	\mathcal S:=\sum_{j,h,k,\ell}
	\gammacompleto
	H_g(x_j,r_{jh}^\supply;x_k,r_{k\ell}^\demand)
\end{equation}
under the constraints $\gammacompleto\geq 0$, $\sum_{jhk\ell}\gammacompleto=1$, and
for the homogeneous marginals
that 
\[
m_j^\supply = \sum_{h}
r_{jh}^\supply \sum_{k\ell} \gammacompleto\,,\qquad
m_k^\demand= \sum_{\ell}
r_{k\ell}^\demand \sum_{jh}\gammacompleto\,,
\]
and
\[
\sum_h r_{jh}^\supply=m_j^\supply\,,\qquad
\sum_\ell r_{k\ell}^\demand =m_k^\demand\,,
\]
for all $j,k\in\{1,\ldots,N\}$ (see also Figure \ref{fig:SSbasics}).
We point out that this latter approach differs from that of the previous one only numerically, because the distance defined in this way in principle coincides with the one discussed in Section \ref{sec:GWD-savare-discrete}.

\subsubsection{Computational approach}\label{sec:GWD-savare-numerica}
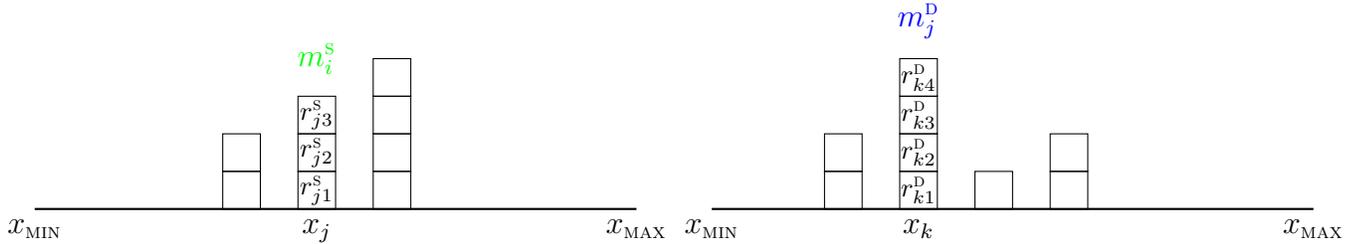
\begin{figure}[h!]
	\centering
	\begin{tikzpicture}[scale=1.0]
	\draw[thick] (-8.5,0) -- (-0.5,0);
	\node at (-8.5,0) [below] {$\xmin$};
	\node at (-0.5,0) [below] {$\xmax$};
	\draw[thick] (0.5,0) -- (8.5,0);
	\node at (0.5,0) [below] {$\xmin$};
	\node at (8.5,0) [below] {$\xmax$};
	\draw (-6,0) rectangle (-5.5,0.5); 
	\draw (-6,0.5) rectangle (-5.5,1);
	\node at (-4.75,0)[below]{$x_i$}; 
	\node at (-4.75,2){\color{cyan} $m_i^\supply$};
	\draw (-5,0) rectangle (-4.5,0.5); 
	\draw (-5,0.5) rectangle (-4.5,1); 
	\draw (-5,1) rectangle (-4.5,1.5);
	\node at (-4.75,0.25){\footnotesize $r^\supply_{i1}$};
	\node at (-4.75,0.75){\footnotesize $r^\supply_{i2}$};
	\node at (-4.75,1.25){\footnotesize $r^\supply_{i3}$};
	\draw (-4,0) rectangle (-3.5,0.5); 
	\draw (-4,0.5) rectangle (-3.5,1); 
	\draw (-4,1) rectangle (-3.5,1.5);
	\draw (-4,1.5) rectangle (-3.5,2);
	\draw (2,0) rectangle (2.5,0.5); 
	\draw (2,0.5) rectangle (2.5,1); 
	\node at (3.25,0)[below]{$x_j$}; 
	\node at (3.25,2.5){\color{orange} $m_j^\demand$};
	\node at (3.25,0.25){\footnotesize $r^\demand_{j1}$};
	\node at (3.25,0.75){\footnotesize $r^\demand_{j2}$};
	\node at (3.25,1.25){\footnotesize $r^\demand_{j3}$};
	\node at (3.25,1.75){\footnotesize $r^\demand_{j4}$};
	\draw (3,0) rectangle (3.5,0.5); 
	\draw (3,0.5) rectangle (3.5,1); 
	\draw (3,1) rectangle (3.5,1.5); 
	\draw (3,1.5) rectangle (3.5,2); 
	\draw (4,0) rectangle (4.5,0.5); 
	\draw (5,0) rectangle (5.5,0.5); 
	\draw (5,0.5) rectangle (5.5,1); 
\end{tikzpicture}
	\caption{Additional discretization of the concentrated masses.}
	\label{fig:SSbasics}
\end{figure}

Starting again from the basic numerical setting introduced in Section \ref{sec:WD-numerica}, we introduce the additional ``vertical'' discretization at any fixed grid point $x_i$ and for each mass distribution, see Figure  \ref{fig:SSbasics}.

At each grid point $x_i$, the mass $m_i^\supply$ is divided in $n^\supply_i\geq 0$ intervals (reasonably, $n^\supply_i=0$ iff $m_i^\supply=0$). 
Analogously, the mass $m_i^\demand$ is divided in $n^\demand_i\geq 0$ intervals (again, $m^\demand_i=0$ iff $m_i^\demand=0$). 
Note that $n^{\{\supply,\demand\}}_i$ acts as additional discretization parameter, and we expect the numerical error to vanish only for both $n^{\{\supply,\demand\}}_i\to +\infty$ and $N\to +\infty$.

Now, at each subdivision of the supply mass is associated the additional variable $\rjh$, $j=1,\ldots,N$, $h=1,\ldots, n^\supply_j$, and at each subdivision of the demand mass is associated the additional variable $\skl$, $k=1,\ldots,N$, $\ell=1,\ldots, n^\demand_k$.

Define also 
$$
\Gamma\in\R^d,\quad \text{with}\quad d=\sum_{i=1}^N n^\supply_i \cdot \sum_{i=1}^N n^\demand_i,
$$
as the vector of all the weights $\gammacompleto$'s associated to the pairs $(\rjh,\skl)$, $j=1,\ldots,N$, $h=1,\ldots, n^\supply_j$, $k=1,\ldots,N$, $\ell=1,\ldots, n^\demand_k$.

The optimization problem to be solved reads as
\begin{equation}\label{OP_SS}
	\begin{tabular}{ll}
		minimize & 
		$\mathcal S(\Gamma,r^\supply,r^\demand)$
		\\ [3mm]
		subject to 
		& $\displaystyle \gammacompleto\geq 0, \quad \forall j,h,k,\ell, \qquad \sum_{jhk\ell}\gammacompleto=1$, \\  [3mm]
		& $\displaystyle m^\supply_j=\sum_{h=1}^{n^\supply_j}\Big(\rjh \sum_{k\ell}\gammacompleto\Big), \quad \forall j$, \\  [4mm]
		& $\displaystyle m^\demand_k=\sum_{\ell=1}^{n^\demand_k}\Big(\skl \sum_{jh}\gammacompleto\Big), \quad \forall k$, 
	\end{tabular}
\end{equation}
where $\mathcal S$ is defined in \eqref{def:S}, and the sums $\sum_{jhk\ell}$, $\sum_{k\ell}$, $\sum_{jh}$ run on all over the index combinations.
It is important to note that the linear systems for $r^\supply$, $r^\demand$ appearing in \eqref{OP_SS} can be actually underdetermined for a given $\Gamma$.
Analogously to the previous cases, we will denote the solution of the optimization problem by $\tilde W^\sss_2$.

\medskip

In order to solve the optimization problem \eqref{OP_SS}, we have devised two strategies: the first one is based on a slow but exhaustive search, which guarantees to find the global minimum. 
The second one is based on a nondifferential descend method and allows for faster computations, although \emph{a priori} limited to the search of a local minimum.

\paragraph{Algorithm SS-A: Exhaustive search}
\begin{enumerate}
	\item Introduce a grid with $Q$ nodes in the unit interval $[0,1]$, namely $0=q_1<q_2<\ldots,q_Q=1$. 
	Note that $Q$ is yet another approximation parameter (the associated error is expected to vanish for $Q\to +\infty$).
	\item List all the possible vectors $\{\tilde\Gamma\in\R^d\}$ such that each element belongs to the set $\{q_1,\ldots,q_Q\}$, and the sum of the elements equals 1. A good algorithm to do this is the following (with $d=5$), although some generated $\tilde\Gamma$'s lead to linear systems with no solution (we have discarded them).
	\\
	\begin{algorithmic}[1]
		\For {$i_1=0,\ldots,Q$}
		\For {$i_2=0,\ldots,Q-i_1$}
		\For {$i_3=0,\ldots,Q-i_1-i_2$}
		\For {$i_4=0,\ldots,Q-i_1-i_2-i_3$}
		\State $i_5=Q-i_1-i_2-i_3-i_4$
		\State $\tilde\Gamma=(i_1,i_2,i_3,i_4,i_5)/Q$
		\EndFor
		\EndFor
		\EndFor
		\EndFor
	\end{algorithmic} 
	%
	%
	\item For each $\tilde\Gamma$, find $\{\rjh\}_{jh}$ and $\{\skl\}_{k\ell}$ (depending on that $\tilde\Gamma$) solving the underdetermined linear systems, selecting the solution with minimal norm.
	\item For each $\tilde\Gamma$ (and corresponding $\{\rjh\}_{jh}$, $\{\skl\}_{k\ell}$), evaluate $\mathcal S(\tilde\Gamma)$.
	\item Find $\min_{\tilde\Gamma}\mathcal S(\tilde\Gamma)$ by exhaustive comparison.
\end{enumerate}

\paragraph{Algorithm SS-B: Random guess and random descend}
\begin{enumerate}
	\item List a largest-as-possible number of vectors $\tilde\Gamma\in\R^d$ such that each element belongs to the interval $[0,1]$ and the sum of the elements equals 1. 
	A simple algorithm to do this is the following.
	\\
	\begin{algorithmic}[1]
		\For {$z=1,\ldots,d$}
		\State Extract a uniform random variable $\omega_z\in\{1,\ldots,Q\}$
		\EndFor
		\State $\tilde\Gamma=(\omega_1,\ldots,\omega_d)/\sum_z \omega_z$
	\end{algorithmic} 
	%
	\item For each $\tilde\Gamma$, find $\{\rjh\}_{jh}$ and $\{\skl\}_{k\ell}$ (depending on that $\tilde\Gamma$) solving the underdetermined linear systems, selecting the solution with minimal norm.
	\item For each $\tilde\Gamma$ (and corresponding $\{\rjh\}_{jh}$, $\{\skl\}_{k\ell}$), evaluate $\mathcal S(\tilde\Gamma)$.
	\item Find $\tilde\Gamma^*:=\arg\min_{\tilde\Gamma}\mathcal S(\tilde\Gamma)$ by exhaustive comparison. 
	\item Starting from $\tilde\Gamma^*$ as initial guess, explore the control space perturbing the best $\Gamma$ obtained so far with small perturbation $\pm \varepsilon$ on random pairs of its elements (thus keeping the sum equal to 1). 
	If the perturbation leads to a lower value, the best $\Gamma^*$ is updated with the perturbed one, otherwise the attempt is discarded. 
	As usual, duly diminishing $\varepsilon$ while the algorithm runs can improve the accuracy (simulated annealing).
\end{enumerate}

In Appendix \ref{sec:GWD-SS-numericaltest}, we present some academic examples to illustrate and better understand the two approaches described in Sects.\ \ref{sec:GHK-numerica} and \ref{sec:GWD-savare-numerica}.

\subsection{Comparison of the three GWDs}\label{sec:comparison}
In this numerical test, we compare the four computational approaches for the three GWDs individually studied above.
This is done by devising a numerical test which can be run by all algorithms.
From now on, the four algorithms will be referred to by FG, PR, SS (variant B only), and GHK, with obvious meaning.

Denoting by $\chi$ the indicator function, we consider the case of two time-dependent distributions with parameters $\alpha,\beta,\eta>0$,
\begin{equation}\label{def:distributions_cfr}
	\rho^\supply(x,t)=\alpha\chi_{\left[-t-\eta,-t+\eta\right]}(x), \qquad
	\rho^\demand(x,t)=\beta\chi_{\left[t-\eta,t+\eta\right]}(x), \qquad 0\leq t\leq T
\end{equation}
in the domain $X=[-T-\eta,T+\eta]$.
The two step functions start with perfectly overlapping support, then supply mass moves leftward while demand mass moves rightward, until they both reach the boundary of the domain, see Figure \ref{fig:test-cfr-spiegazione}.
\begin{figure}[h!]
	\centering
	\begin{tikzpicture}[scale=1.0]
	\draw[thick] (-8.5,0) -- (-0.5,0);
	\node at (-8.1,0) [above] {\tiny supply};
	\draw[thick] (0.5,0) -- (8.5,0);
	\node at (8.0,0) [above] {\tiny demand};

	\draw[thick] (-5,0) rectangle (-3,1); 
	\node at (-5,0) [below] {\footnotesize $-t-\eta$}; 
	\node at (-3,0) [below] {\footnotesize $-t+\eta$}; 
	\node at (-5,1) [left] {$\alpha$}; 
	\draw[->] (-5.1,0.5) -- (-6,0.5); 

	\draw[thick] (3,0) rectangle (5,1.5); 
	\node at (3,0) [below] {\footnotesize $t-\eta$}; 
	\node at (5,0) [below] {\footnotesize $t+\eta$}; 
	\node at (3,1.5) [left] {$\beta$}; 
	\draw[->] (5.1,0.75) -- (6,0.75); 
\end{tikzpicture}
	\caption{Moving distributions $\rho^\supply$ and $\rho^\demand$ as defined in \eqref{def:distributions_cfr}.}
	\label{fig:test-cfr-spiegazione}
\end{figure}

In particular, for this comparison test we set $\alpha\in\{1,2\}$, $\beta=2$, $\eta=1$, $T=4$. Moreover, we have the following algorithm-specific parameters:
\begin{itemize}
	\item PR: $a=2.5$, $b=1$.
	\item GHK: $a=b=1$. 
	\item SS: $a=b=1$, 4 atomic masses per distribution, 4 additional divisions per mass, $Q=1,\!000$, $\epsilon=0.001$, iterations 5,000 (random searches) + 40,000 (random descend), average on 10 runs.
\end{itemize}

\begin{figure}[h!]
	\centering
	\includegraphics[width=0.49\textwidth]{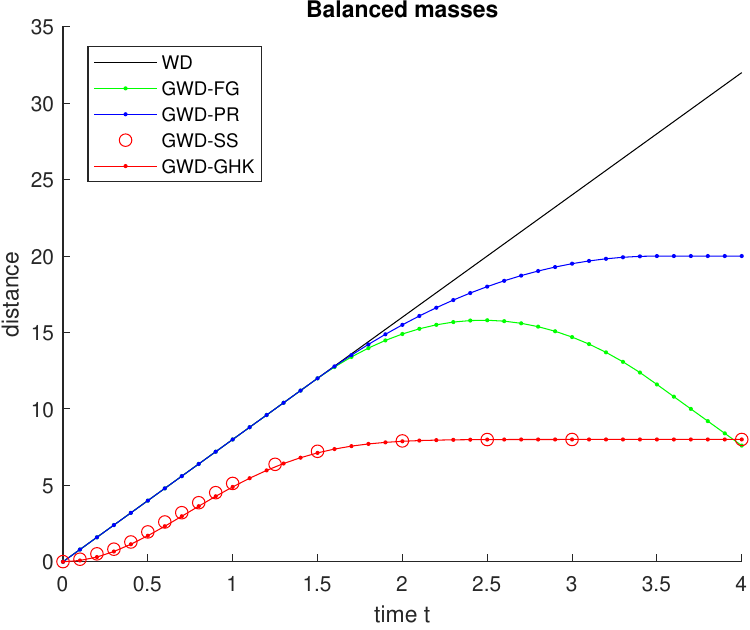}
	\includegraphics[width=0.49\textwidth]{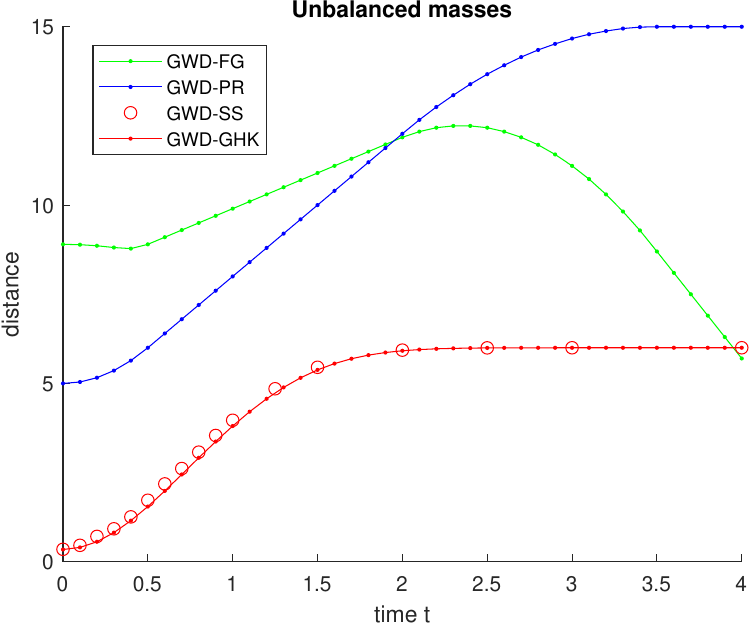}
	\caption{
		Comparison of the three GWDs between $\rho^\supply$ and $\rho^\demand$ as defined in \eqref{def:distributions_cfr}, as a function of time $t$ ($\beta=2$, $\eta=1$, $T=4$). 
		$\alpha=2$ (balanced), including also classical WD (left), and $\alpha=1$ (unbalanced, right). 
	}
	\label{fig:test-cfr}
\end{figure}

Figure \ref{fig:test-cfr} shows the behavior of the three GWDs between $\rho^\supply$ and $\rho^\demand$ as defined in \eqref{def:distributions_cfr}, as a function of time. 
WD is also shown in the balanced case $\alpha=\beta$.
%

This comparison test allows us to sketch some preliminary conclusions.
First of all, SS seems to be computationally unfeasible. It is impossible to deal with a reasonable number of concentrated masses and additional divisions, therefore it will be discarded in the following analysis.   

Among the others, PR and GHK share with $L^p$ an important drawback: the distances saturate when distributions do not overlap and they are far enough from each other. This does not seem to be a nice feature in the context of traffic flow models.

Lastly, FG has an important feature which deserves a discussion: the distance decreases whenever distributions are close to the boundaries, even if they are very far from each other. On first glance, this seems to be an issue, but it can be also be seen as a good point: basically, it implies that what happens near the boundaries is less important than what happens on the rest of the road. 
Considering the fact that boundary conditions are typically unknown in traffic flow modeling, it could make sense to give priority to that part of the road which is less or not at all influenced by boundary conditions.  

%
%
%
%
%

\section{Applications to traffic model analysis}\label{sec:traffic}
In this section we will try to understand to what extent the GWDs previously studied are suitable in the context of traffic flow models.
It is useful to recall here that in \cite{briani2018CMS} authors have already extensively studied the sensitivity of the LWR model \cite{lighthill1955PRSLA, richards1956OR} using the WD. 
More precisely, many simulations were performed varying some key parameters like the initial vehicle density, the road capacity (maximal flux), the critical density (vehicle density associated to the maximal flux), the size of the road network, and the distribution coefficients at junctions (how traffic distributes across a junction with more than one outgoing road).
Interestingly, authors found that the parameter the model is most sensitive to is the last one, namely the distribution coefficients at junctions.
That said, the usage of WD have limited the authors to the comparison of \emph{balanced} vehicle distributions only. 
This means, in practice, that the number of vehicles should be  the same in the two considered scenarios. 

In the following we try to complete the sensitivity analysis including the cases in which the mass distributions are unbalanced.

\subsection{Mathematical modeling}\label{sec:trafficmodelling}
In a macroscopic setting, traffic flow is described by means of average quantities only, like \emph{density} $\rho$ (vehicles per unit of space, at a given time), \emph{flux} $f$ (vehicles per unit of time crossing a given point, at a given time), and \emph{velocity} $v$.
Models can be \emph{first-order} (i.e., velocity based), if velocity is given as a function of the density and acceleration is neglected, or \emph{second-order} (i.e., acceleration based) if, instead, velocity is an independent variable and acceleration is given as a function of both density and velocity.

In the following we will assume the road is a finite segment $X=[\xmin,\xmax]$ with a single lane and a single class of vehicles. 

\subsubsection*{LWR model}
The LWR model, introduced independently by Lighthill \& Whitham \cite{lighthill1955PRSLA} and by Richards \cite{richards1956OR}, describes the evolution of the vehicle density $\rho(x,t)$ by means of the following IBVP involving an hyberbolic PDE
\begin{equation}\label{LWR}
	\left\{
	\begin{array}{ll}
		\partial_t \rho + \partial_x f(\rho)=0, & t\in[t_0,T], \quad x\in [\xmin,\xmax], \\
		\rho(x,0)=\rho_0(x), & x\in [\xmin,\xmax], \\
		\rho(\xmin,t)= \rhoin(t), & t\in[t_0,T], \\
		\rho(\xmax,t)= \rhoout(t), & t\in[t_0,T], 
	\end{array}
	\right.
\end{equation}
for some initial time $t_0$, final time $T>t_0$, and concave flux function $f:[0,\rhomax]\to\R^+$ (called \emph{fundamental diagram} in the context of traffic flow) such that $f(0)=f(\rhomax)=0$, for some maximal density $\rhomax$. 
Note that we have naturally guaranteed the constraint $0\leq \rho\leq \rhomax$ at any time, provided $0\leq \rho_0,\rhoin,\rhoout\leq \rhomax$.

Recalling the well known physical law $f=\rho v$, which is always valid for any $x$ and $t$, we get that the LWR model implicitly assumes that the vehicle velocity $v$ only depends on the density $\rho$ (since $f=f(\rho)$, therefore $v=v(\rho)$).
From the modeling point of view, this means that the velocity adapts instantaneously to the changes of the density, i.e.\  acceleration is infinite. For this reason the LWR model is said to be \emph{first-order}.

As usual in the mathematical literature, in the following Tests 1, 2, 3 we will define 
\begin{equation}\label{def:DF}
	f(\rho)=\rho(1-\rho),
\end{equation}
implicitly assuming that the maximal density and the maximal velocity of vehicles are normalized and equal to 1. 
This is done for simplicity, since a qualitative analysis is sufficient to highlight the differences among the GWDs.
In Test 4, instead, we consider real-like parameters.

\medskip

Equation \eqref{LWR} can be numerically approximated by means of, e.g., the Godunov scheme.
Using again the numerical setting introduced in Section \ref{sec:WD-numerica}, this scheme reads as
\begin{equation}\label{Godunovscheme}
	\rho_i^{n+1}=\rho_i^n-\frac{\Dx}{\Dt}\Big(\mathcal G(\rho_i^n,\rho_{i+1}^n)-\mathcal G(\rho_{i-1}^n,\rho_i^n) \Big),\qquad i=2,\ldots,N-1, \quad n=0,\ldots, \left[\frac{T-t_0}{\Dt}\right],
\end{equation}
where $\rho_i^n$ represents the approximated value of $\rho$ at the center of the cell $C_i$ at time $t_0+n\Dt$, and $\mathcal G$ is the Godunov's numerical flux defined by
\begin{equation}\label{Godunovflux}
	\mathcal G(\rhomeno,\rhopiu):=
	\left\{
	\begin{array}{ll}
		\displaystyle \min_{z\in [\rhomeno,\rhopiu]} f(z), & \textrm{if } \rhomeno\leq \rhopiu, \\ [3mm]
		\displaystyle \max_{z\in [\rhopiu,\rhomeno]} f(z), & \textrm{if } \rhomeno\geq \rhopiu. \\
	\end{array}
	\right.
\end{equation}
\begin{rmk}\label{rem:fluxes}
	Boundary conditions can be also given in form of incoming and outgoing flux, rather than as left and right densities $\rhoin$, $\rhoout$. 
	This happens, e.g., if boundary data are given by fixed sensors counting vehicles passing through $x=\xmin$ and $x=\xmax$ in a time interval $\Dt$.
	If this is the case, values $\rho_1^n=\rhoin^n$ and $\rho_N^n=\rhoout^n$ are dropped, and incoming/outgoing fluxes $\fin^n$, $\fout^n$ are directly injected in the scheme \eqref{Godunovscheme} taking the place of 
	$\mathcal G(\rho_{i-1}^n,\rho_{i}^n)$ if $i=1$, 
	and of 
	$\mathcal G(\rho_i^n,\rho_{i+1}^n)$ if $i=N$, respectively.
\end{rmk}

\subsubsection*{ARZ model}
The ARZ model, introduced independently by Aw \& Rascle \cite{aw2000SIAP} and by Zhang \cite{zhang2002TRB}, describes the evolution of the vehicle density $\rho(x,t)$ and the velocity $v(x,t)$ by means of the following IBVP involving two hyberbolic PDEs
\begin{equation}\label{ARZ}
	\left\{
	\begin{array}{ll}
		\partial_t \rho + \partial_x (\rho v)=0, & t\in[t_0,T], \quad x\in [\xmin,\xmax], \\
		\partial_t (\rho w) + \partial_x (v\rho w)=\rho \frac{v_{eq}(\rho)-v}{\tau}, & t\in[t_0,T], \quad x\in [\xmin,\xmax], \\
		\rho(x,0)=\rho_0(x), & x\in [\xmin,\xmax], \\
		v(x,0)=v_0(x), & x\in [\xmin,\xmax], \\
		\rho(\xmin,t)= \rhoin(t), & t\in[t_0,T], \\
		\rho(\xmax,t)= \rhoout(t), & t\in[t_0,T], \\
		v(\xmin,t)= \vin(t), & t\in[t_0,T], \\
		v(\xmax,t)= \vout(t), & t\in[t_0,T],
	\end{array}
	\right.
\end{equation}
where $w(x,t):=v(x,t)+P(\rho(x,t))$, 
$
P(\rho):=\frac{v_{ref}}{\gamma}\left(\frac{\rho}{\rhomax}\right)^\gamma,
$
and $v_{ref}$, $\gamma$, an $\tau$ are strictly positive parameters.
This time, the velocity $v$ is an independent unknown and the acceleration $\frac{v_{eq}(\rho)-v}{\tau}$ appears explicitly in the model, depending on $\rho$ and $v$. 
For this reason, the ARZ model is said to be \emph{second-order}.

In order to fairly compare the two models we define $v_{eq}(\rho):=\frac{f(\rho)}{\rho}$, so that the desired equilibrium velocity of the ARZ model coincides with the instantaneous velocity of the LWR model. Also in the boundary conditions, we set $\vin=v_{eq}(\rhoin)$ and $\vout=v_{eq}(\rhoout)$.

Equation \eqref{ARZ} is numerically approximated by means of the Lax--Friedrichs scheme.

%
%
%
%
%
%
\subsection{Numerical tests}

\subsubsection{Test 1: sensitivity to boundary conditions}\label{Test1traffic}
In this test we investigate the sensitivity to boundary conditions.
We consider the case described in Remark \ref{rem:fluxes}, which is a scenario often met in real life \cite{brianiAVpp}.
In the same spirit of \cite{brianiAVpp}, we assume that the time interval is divided in two subintervals: $[t_0,0]$ and $(0,T]$, where $t=0$ corresponds to the current time (the `now'), i.e.\ the time the simulation is actually performed. 
Until time 0, we can rely on the sensors' measurements $\fin$, $\fout$, and we use it as boundary conditions (see Remark \ref{rem:fluxes}). 
After that time, we have no such a data because they are not yet available, but we assume to be able to (exactly) forecast the average value of them,
$$
\pfin=\frac1T \int_0^T \fin(t)dt,\qquad
\pfout=\frac1T \int_0^T \fout(t)dt.
$$
In short: we denote by $\rho^\supply$ the ``exact'' solution obtained by using boundary fluxes $\fin(t)$, $\fout(t)$ in the whole interval $[t_0,T]$, and by $\rho^\demand$ the ``predicted'' solution obtained by using flux boundary data $\fin(t)$, $\fout(t)$ until time $t=0$ and then by using the predicted average values $\pfin$, $\pfout$ until time $T$.
At any time $t$, the vehicle distributions $\rho^\supply(\cdot,t)$ and $\rho^\demand(\cdot,t)$ are compared computing the GWDs. 

Initial conditions are $\rho^\supply(x,t_0)=\rho^\demand(x,t_0)=0.2$ for any $x$.
Boundary fluxes are randomly generated: 
$\fin$ has uniform distribution in $[0,\frac14]$, 
$\fout$ has uniform distribution in $[0,\mathcal G(0.2,0.2)]$.
Other parameters are: $\xmin=0$, $\xmax=40$, $t_0=-10$, $T=30$, $\Dx=0.4$ (corresponding to $N=100$), and $\Dt=0.25$.

\begin{figure}[h!]
	\centering
	\includegraphics[width=0.49\textwidth]{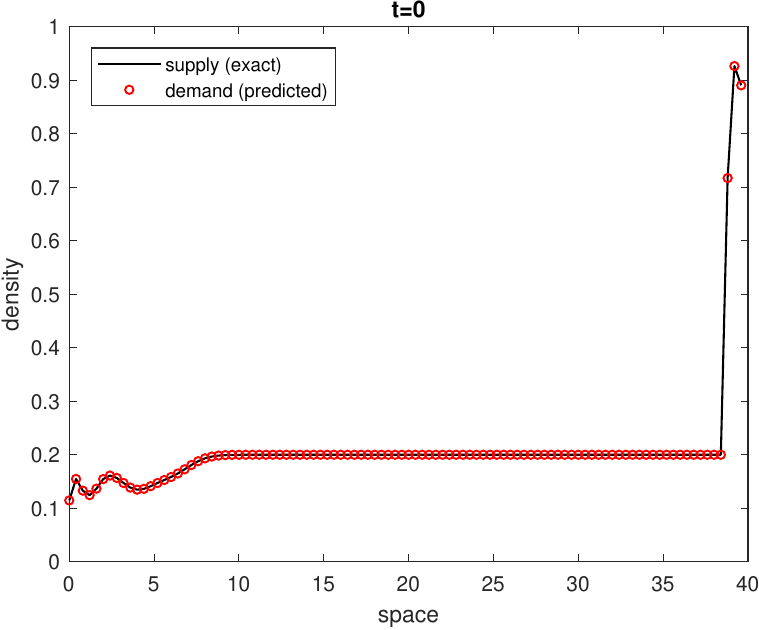}
	\includegraphics[width=0.49\textwidth]{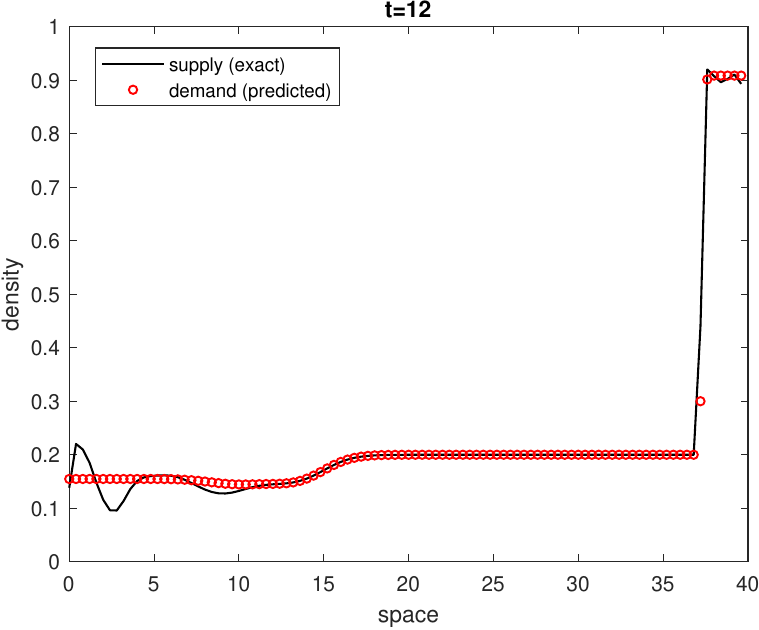} \\
	\caption{
		Test 1 on sensitivity to boundary conditions. Two screenshots at $t=0$ (`now'), i.e.\ when solutions cease to be equal (left), and $t=12$ (right).
	}
	\label{fig:trafficT1-screenshots}
\end{figure}

\begin{figure}[h!]
	\centering
	\includegraphics[width=0.49\textwidth]{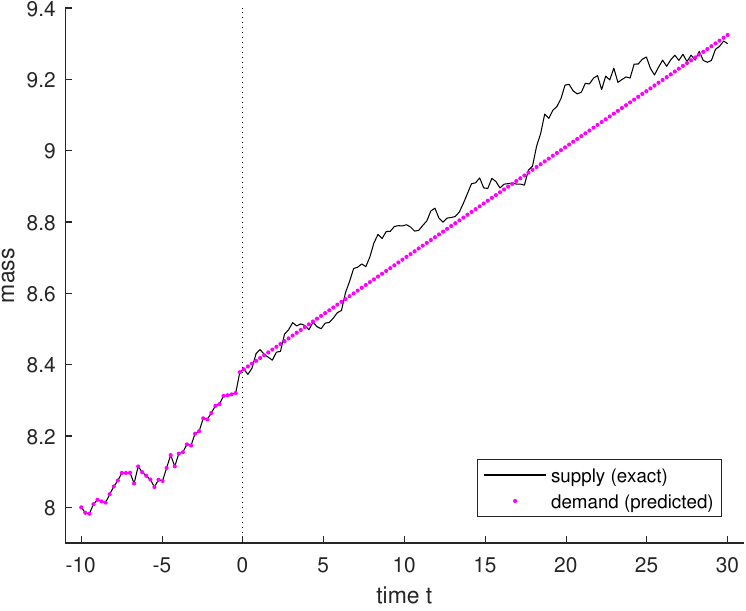}
	\includegraphics[width=0.49\textwidth]{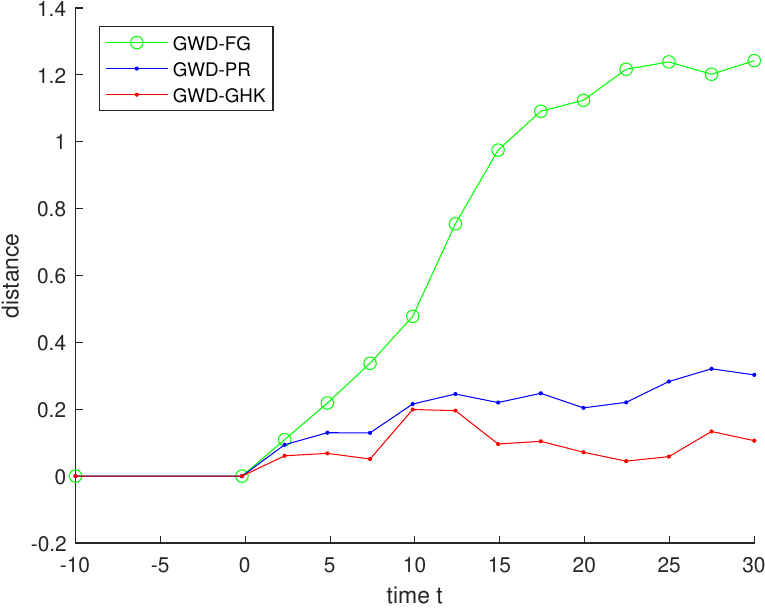}
	\caption{
		Test 1 on sensitivity to boundary conditions. Difference of mass between the two simulations (left) and comparison of GWDs (right).
	}
	\label{fig:trafficT1-masses-and-GWDs}
\end{figure}

Figure \ref{fig:trafficT1-screenshots} shows two screenshots of the simulations at times $t=0$ and $t=12$,
while 
Figure \ref{fig:trafficT1-masses-and-GWDs} shows the total mass on the road and the GWDs between $\rho^\supply(\cdot,t)$ and $\rho^\demand(\cdot,t)$ at any time.

Observing the numerical results, we tend to confirm what we have anticipated in Section \ref{sec:comparison}, namely that the Figalli \& Gigli approach seems to be the most suitable GWD to be used in the context of traffic flow modeling, mainly because of how boundaries are regarded.
In this test, the discrepancy between the exact and predicted solution slowly propagates from the boundaries to the interior of the road, impacting more and more on the  reliability of the traffic prediction. 
This behavior of the error is clearly better caught by FG (see Figure \ref{fig:trafficT1-masses-and-GWDs}(right)), which is the last distance that saturates.
%
%
%
%
\subsubsection{Test 2: sensitivity to traffic light cycle}\label{Test2traffic}
In this test we consider a shorter road with respect to the previous test, and we assume there is a traffic light in the middle of the road, at the interface between two numerical cells. 
The duration of the red phase equals that of the green phase.
Technically, the red traffic light is modeled by imposing null flux at the right (resp., left) boundary of the cell before (resp., after) the traffic light.
The aim is to investigate the sensitivity to the duration of the traffic light cycle. 

Initial conditions are $\rho^\supply(x,t_0)=\rho^\demand(x,t_0)=0$ for any $x$.
Left Dirichlet boundary conditions are $\rhoin^\supply=\rhoin^\demand=0.4$.
Right Dirichlet boundary conditions are $\rhoout^\supply=\rhoout^\demand=0$. 
Other parameters are: $\xmin=0$, $\xmax=4$, $t_0=0$, $T=20$, $\Dx=0.04$ (corresponding to $N=100$), and $\Dt=0.025$.
Both red and green phases last 50 time steps for supply distribution and 40 time steps for demand distribution.

\begin{figure}[h!]
	\centering
	\includegraphics[width=0.49\textwidth]{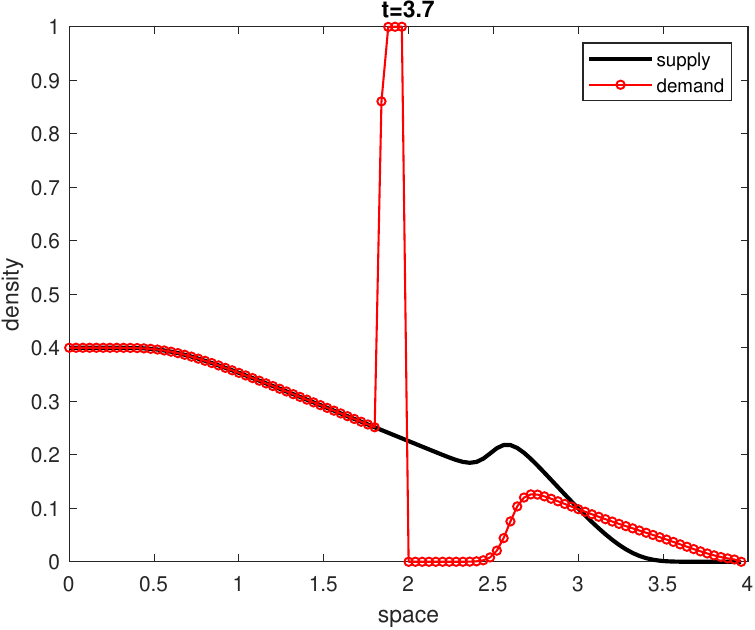}
	\includegraphics[width=0.49\textwidth]{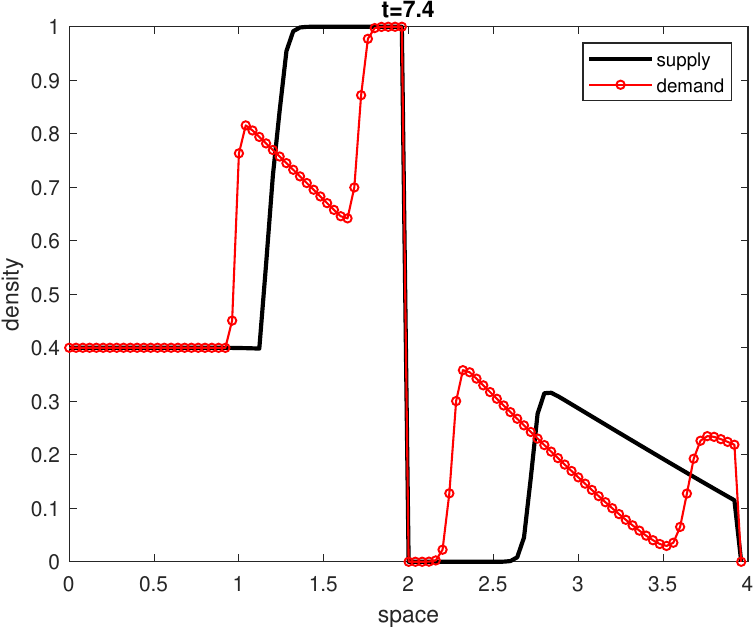} \\
	\caption{
		Test 2 on sensitivity to traffic light cycle. Two screenshots at $t=3.7$ with red traffic light for demand and green traffic light for supply (left), and $t=7.4$ with red traffic light for both distributions (right).
	}
	\label{fig:trafficT2-screenshots}
\end{figure}

Figure \ref{fig:trafficT2-screenshots} shows two screenshots of the simulations at times $t=3.7$ and $t=7.4$.

\begin{figure}[h!]
	\centering
	\includegraphics[width=0.49\textwidth]{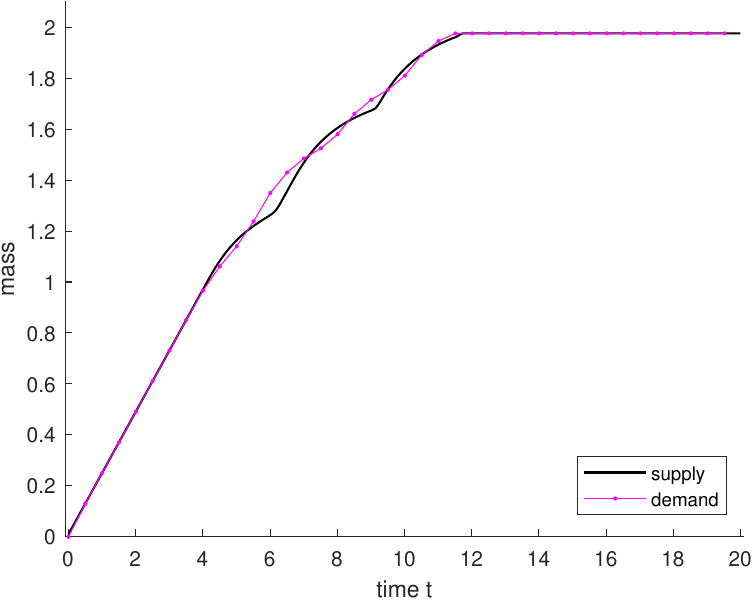}
	\includegraphics[width=0.49\textwidth]{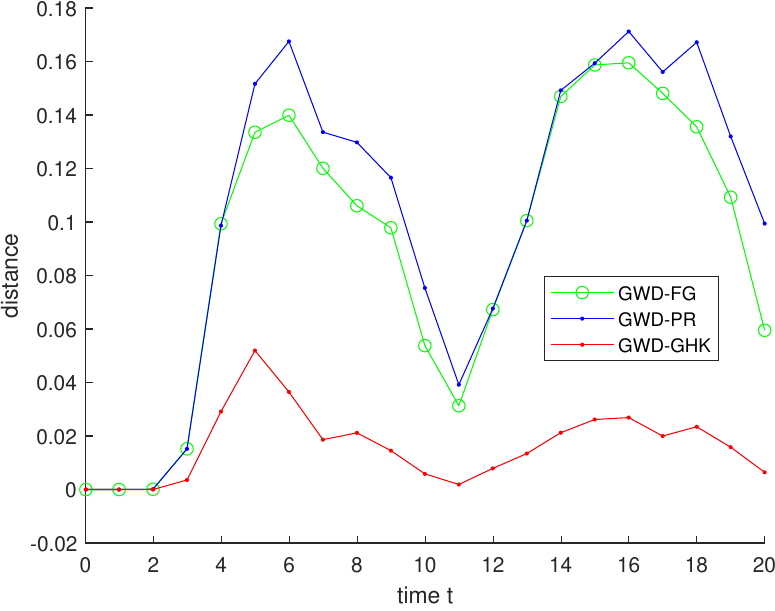}
	\caption{
		Test 2 on sensitivity to traffic light cycle. Difference of mass between the two simulations (correctly masses differ only in a certain time interval, left), and comparison of GWDs (right).
	}
	\label{fig:trafficT2-masses-and-GWDs}
\end{figure}

Figure \ref{fig:trafficT2-masses-and-GWDs} shows the total mass on the road in the two scenarios and the GWDs between $\rho^\supply(\cdot,t)$ and $\rho^\demand(\cdot,t)$ at any time.

This test does not show much difference among the three GWDs, having them basically the same cycling behavior. 
That behavior is correct since the two solutions cyclically move away from each other and then come closer again.
%
%
%
%
%
%
\subsubsection{Test 3: sensitivity to model order}\label{Test3traffic}
In this test we quantify the difference between two simulations obtained by means of the LWR model (supply, first-order) and the ARZ model (demand, second-order).

Initial conditions are $\rho^\supply(x,t_0)=\rho^\demand(x,t_0)=0.7$ for $x\in[1.6,2.4]$ and 0.1 elsewhere.
Left and right Dirichlet boundary conditions are null. 
Other parameters are: $\tau=0.05$, $\gamma=2$, $v_{ref}=1$, $\xmin=0$, $\xmax=4$, $t_0=0$, $T=3$, $\Dx=0.02$ (corresponding to $N=200$), and $\Dt=0.01$.

\begin{figure}[h!]
	\centering
	\includegraphics[width=0.49\textwidth]{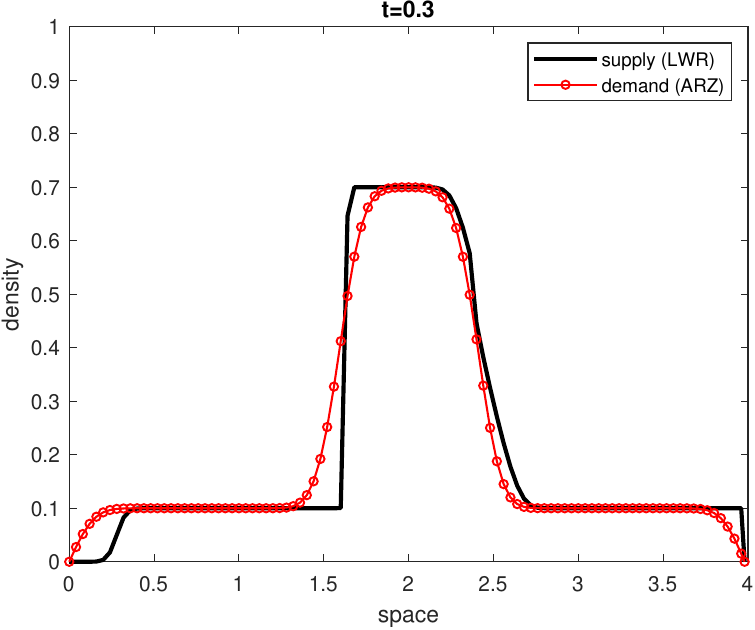}
	\includegraphics[width=0.49\textwidth]{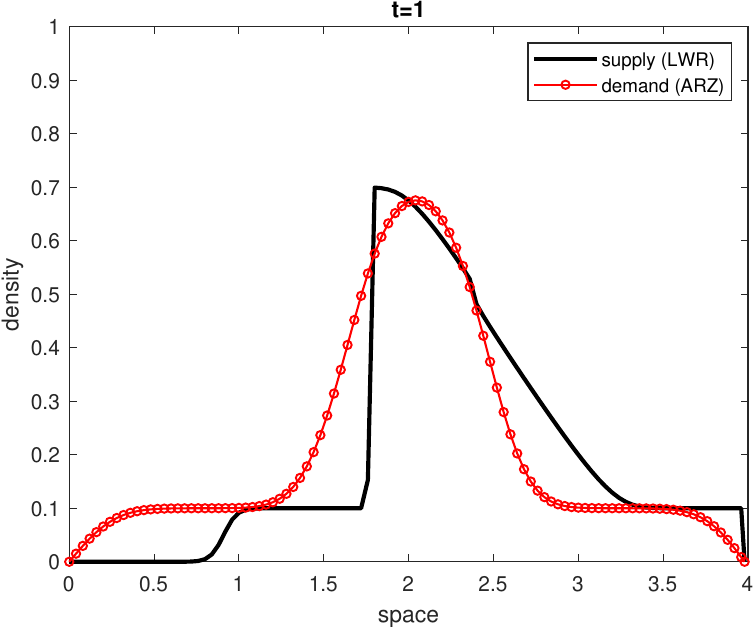} \\
	\caption{
		Test 3 on sensitivity to model order. Two screenshots at $t=0.3$ (left), and $t=1$ (right).
	}
	\label{fig:trafficT3-screenshots}
\end{figure}

Figure \ref{fig:trafficT3-screenshots} shows two screenshots of the simulations at times $t=0.3$ and $t=1$.

\begin{figure}[h!]
	\centering
	\includegraphics[width=0.49\textwidth]{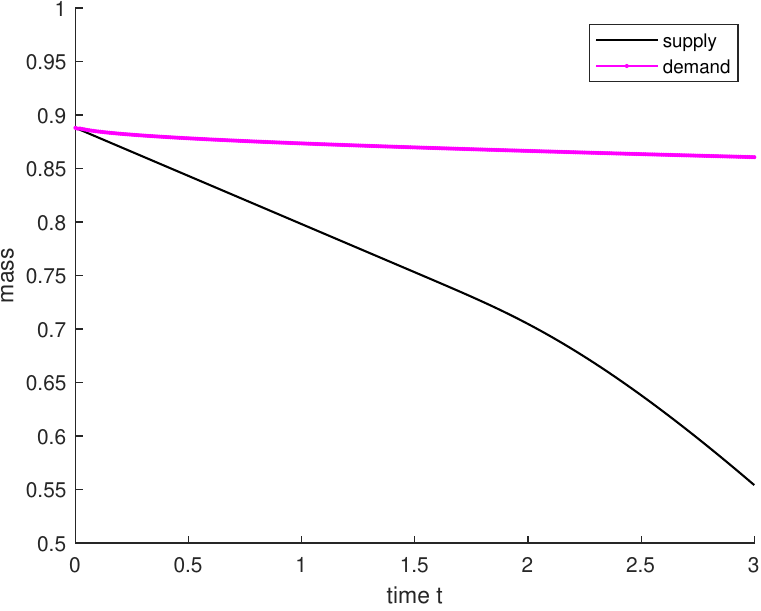}
	\includegraphics[width=0.49\textwidth]{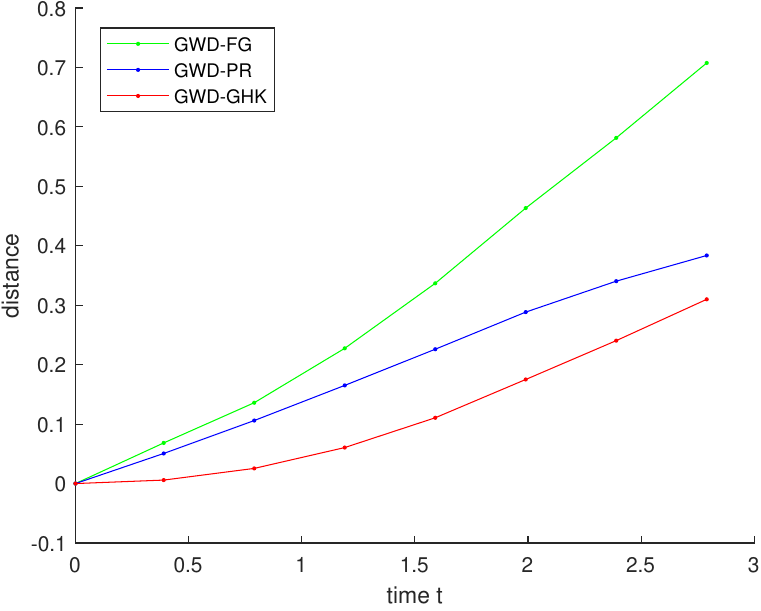}
	\caption{
		Test 3 on sensitivity to model order. Difference of mass between the two simulations (left) and comparison of GWDs (right).
	}
	\label{fig:trafficT3-masses-and-GWDs}
\end{figure}

Figure \ref{fig:trafficT3-masses-and-GWDs} shows the total mass on the road in the two scenarios and the GWDs between $\rho^\supply(\cdot,t)$ and $\rho^\demand(\cdot,t)$ at any time. 
Here we have chosen $a=0.5$ and $b=1$ as parameters for computing $\tilde W^\pr_1$.

This Test is similar to Test 1, even if conclusions are less evident: 
here the two solutions start differing immediately, both at the boundary and inside the road. 
Moreover, intuition tells that the difference increases as the time goes by.
Among the three GWDs, FG is the one which increases more rapidly at the initial time, and it is the one with the major relative increment at later times.

\subsubsection{Test 4: calibration and validation}\label{Test4traffic}
In this section we consider two important problems in traffic modeling, namely \textit{calibration} and \textit{validation}.
In fact, once a new mathematical model is devised, the calibration is needed for finding the best values of its parameters to be used in order to fit as much as possible some reference solution, typically coming from real observations (real traffic data). 
Once the calibration is done, one has to \textit{validate} the model, checking if it is able to compute accurate solutions also in scenarios which are different from those used for calibration.

In both steps the choice of a distance function is crucial, since it is used to measure the difference between the exact and the simulated traffic conditions. 
Different distances can highlight different features, thus affecting the values of the optimal parameters.

To study these aspects of traffic modeling, we resort again to the LWR model (Section \ref{sec:trafficmodelling}). 
Here, we do not normalize maximal density $\rhomax$ and maximal velocity $\velmax$ to 1, so we have  
$$
f(\rho)=\rho\left(-\frac{\velmax}{\rhomax}\rho+\velmax\right)
$$
instead of \eqref{def:DF}.
Summarizing, to run the model we need the two parameters $\rhomax$ and $\velmax$, the initial condition $\rho_0$, and the boundary conditions $\rhoin, \rhoout$.
The proposed numerical tests are performed in two steps:
\begin{itemize}
	\item[-] \emph{direct problem}: we fix the parameters/conditions and we solve \eqref{LWR}, saving the final solution  $\rho_T^{\textsc{ref}}(\cdot):=\rho(\cdot,T)$ as reference solution; 
	\item[-] \emph{inverse problem}: we assume that some parameters/conditions are unknown, and we use the other parameters/conditions and the reference solution $\rho_T^{\textsc{ref}}$ to find the unknown parameters/conditions by means of an optimization problem.  
\end{itemize}
Note that the second step requires to measure several times the distance between the reference solution and the solution obtained with tentative parameters/conditions, until the minimum distance (best match) is found.

\begin{figure}[h!]
	\centering
	\includegraphics[width=0.49\textwidth]{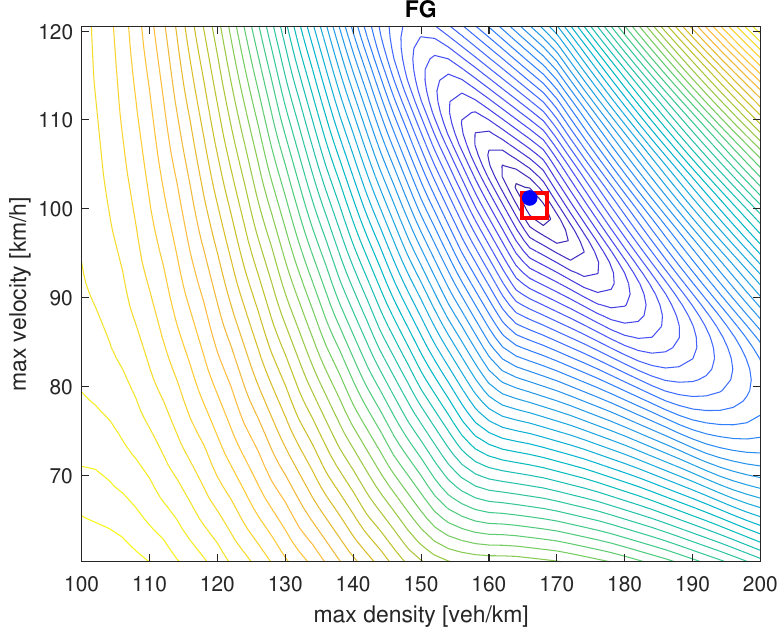}
	\includegraphics[width=0.49\textwidth]{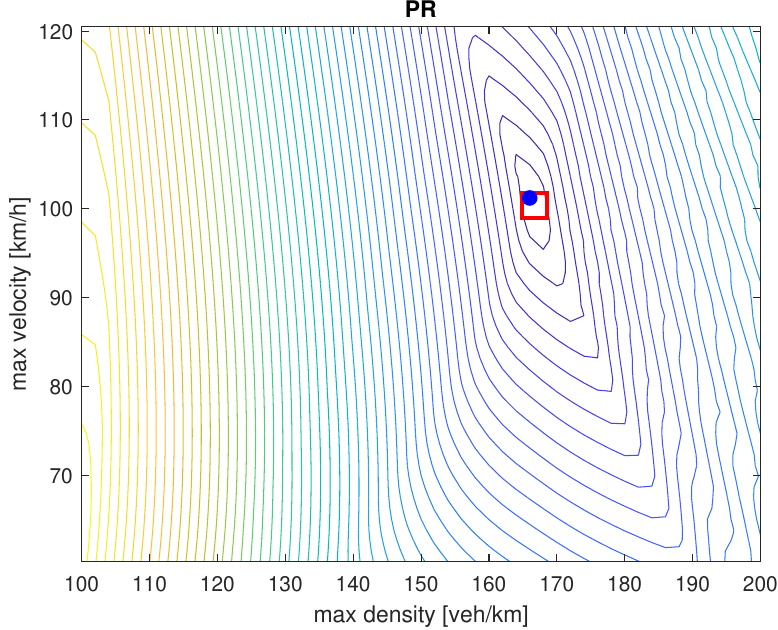}\\
	\includegraphics[width=0.49\textwidth]{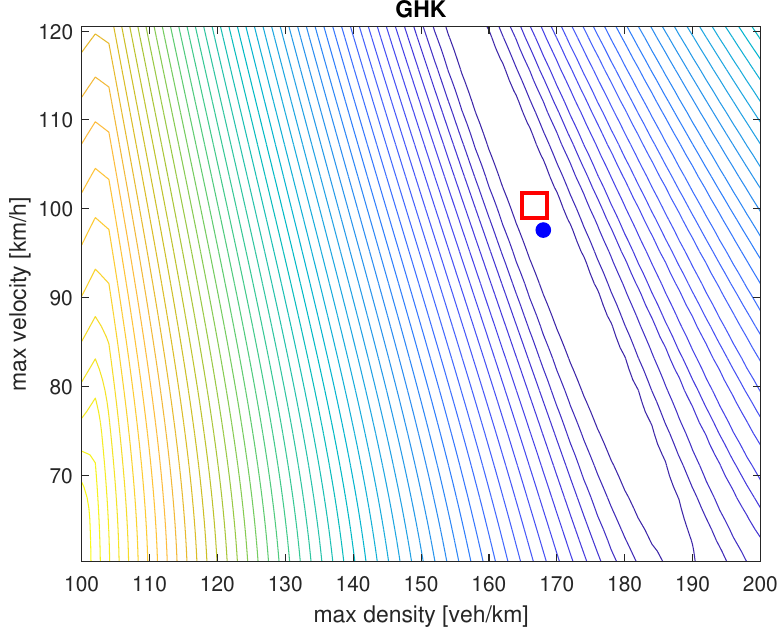}
	\includegraphics[width=0.49\textwidth]{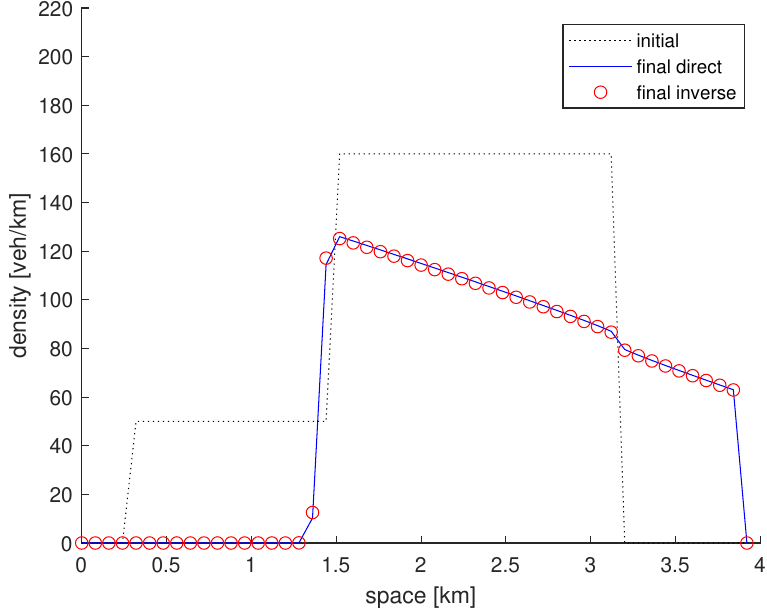} \\
	\caption{
		Test 4a on calibration. Level sets of the GWD functions $(\rhomax,\velmax)\to\tilde W(\rho(\cdot,T;\rhomax,\velmax),\rho_T^{\textsc{ref}})$ for FG (top-left), PR (top-right), and GHK (bottom-left). 
		The red square locates the exact pair $(\rhomax^*,\velmax^*)$, while the blue dot locates the argmin of the function. 
		In the bottom-right figure we show the initial condition $\rho_0$, the reference density function $\rho_T^{\textsc{ref}}$ computed in the direct problem, and the optimal density function $\rho^*$ at final time corresponding to the solution to the optimization problem. 
	}
	\label{fig:trafficT4a}
\end{figure}

\paragraph{Test 4a.}
In this test, for the direct problem we set $\rhomax^*=167$ veh/km, $\velmax^*=100$ km/h, $\xmin=0$ km, $\xmax=4$ km, $T=2$ min, $\rhoin=0$ veh/km, and $\rhoout=0$ veh/km. 
$\rho_0$ is plotted in Figure \ref{fig:trafficT4a}(bottom-right).

For the inverse problem, instead, we assume that $\rhomax$ and $\velmax$ are unknown, and we run an exhaustive search for the best match. 
The intervals for the search are $[100,200]$ veh/km for the density and $[60,120]$ km/h for the velocity.

In Figure \ref{fig:trafficT4a} we show the level sets of the function $(\rhomax,\velmax)\to\tilde W(\rho(\cdot,T;\rhomax,\velmax),\rho_T^{\textsc{ref}})$ for FG, PR, and GHK approaches.

It is interesting to note that in all the three cases the optimization problem is able to find quite precisely the optimal pair $(\rhomax^*,\velmax^*)$, which is indeed the argmin of the GWDs. Nevertheless, the qualitative behavior of the level sets around the minimum point is different. 
FG and PR show the best behavior, with the minimum point well contoured by level sets, meaning that the minimum is well localized. 
Conversely, GHK shows a plateau around the argmin, which clearly makes the minimum less identifiable.

\paragraph{Test 4b}
In this test, for the direct problem we set $\rhomax=167$ veh/km, $\velmax=100$ km/h, $\xmin=0$ km, $\xmax=2$ km, $T=3$ min, $\rhoin^*=40$ veh/km, $\rhoout^*=0$ veh/km, and $\rho_0\equiv 0$. 
Moreover, we assume that between $t=1$ min and $t=1.7$ min an accident at 1.6 km blocks the traffic. 

For the inverse problem, instead, we assume that $\rhoin$ and $\rhoout$ are unknown, and the accident is unknown too (this is a typical situation which can happen during model validation). 
The intervals for the search are $[0,\rhomax/2]$ veh/km for the inflow boundary conditions and $[0,\rhomax]$ km/h for the outflow boundary conditions.

Conversely to the Test 4a, here it is very possible that the optimization problem is not able to find the perfect match with the reference solution since the accident is not at all considered in the inverse problem. Indeed, the optimization problem could tend to compensate the internal discrepancy due to the accident using different boundary conditions. 

\begin{figure}[h!]
	\centering
	\includegraphics[width=0.49\textwidth]{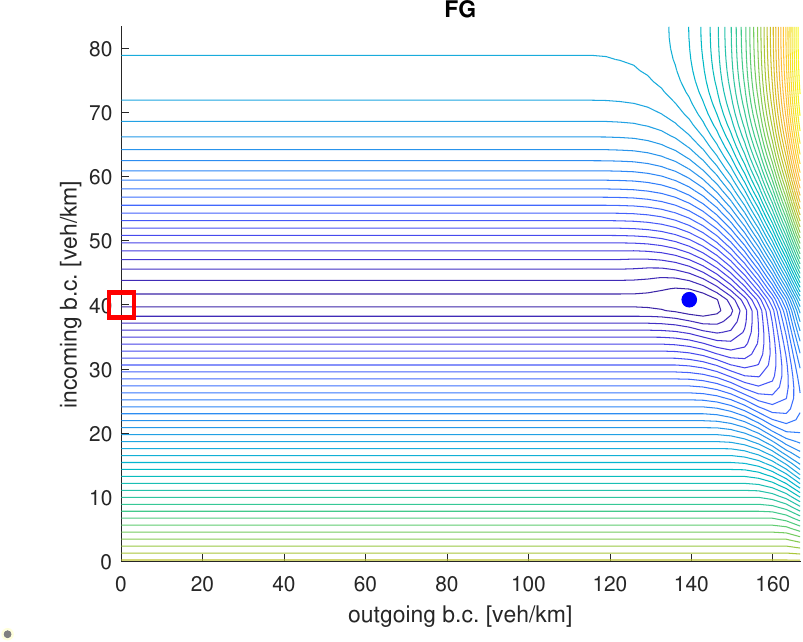}
	\includegraphics[width=0.49\textwidth]{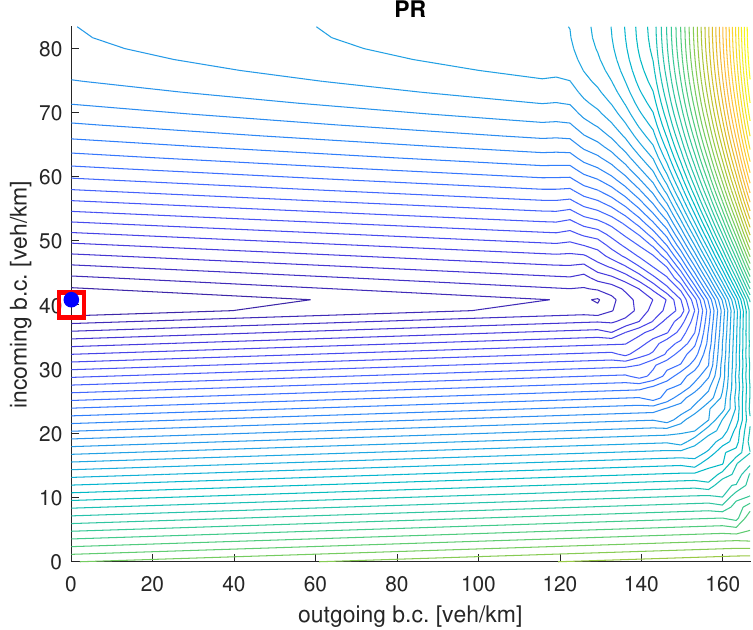}\\
	\includegraphics[width=0.49\textwidth]{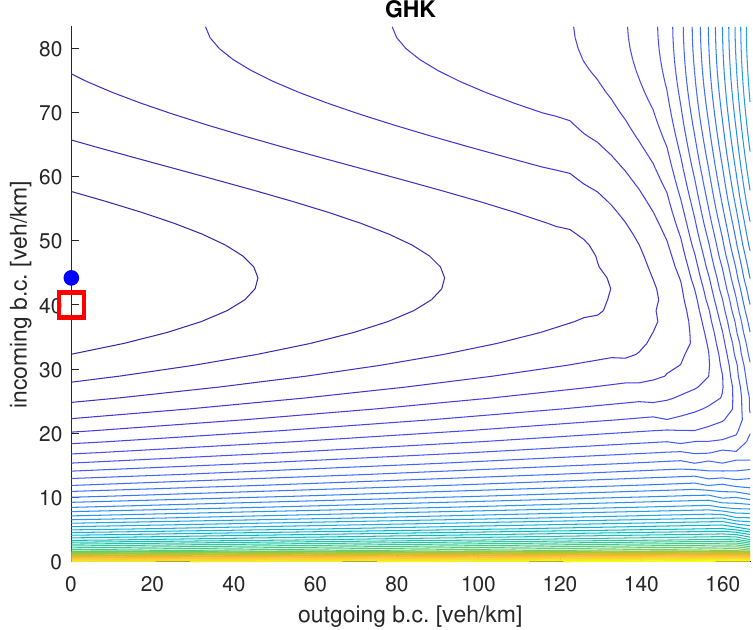}
	\caption{
		Test 4b on validation. Level sets of the GWD functions $(\rhoout,\rhoin)\to\tilde W(\rho(\cdot,T;\rhoin,\rhoout),\rho_T^{\textsc{ref}})$ for FG (top-left), PR (top-right), and GHK (bottom). 
		The red square locates the exact pair $(\rhoout^*,\rhoin^*)$, while the blue dot locates the argmin of the function. 
	}
	\label{fig:trafficT4b}
\end{figure}

In Figure \ref{fig:trafficT4b} we show the level sets of the function $(\rhoout,\rhoin)\to\tilde W(\rho(\cdot,T;\rhoin,\rhoout),\rho_T^{\textsc{ref}})$ for FG, PR, and GHK approaches.

Since the accident happens close to the boundary, and control parameters $\rhoin,\rhoout$ act at boundaries as well, FG struggles to find the optimal parameters, being basically insensitive to any possible variation. This is evident observing the plateau around the minimum, which also lead to a completely wrong estimation of the boundary conditions. 
The other two methods, instead, are more sensitive to variations near boundaries and are able to locate the exact solution even if the accident is ignored.

\section*{Conclusions}
In this paper we have proposed four algorithms for computing, in a discrete context, four GWDs (two of them theoretically equivalent to each other). 
One of them (SS) results to be numerically unfeasible, while the other three (FG, PR, GHK) are all relative simple to implement once the computational procedure is arranged.

Regarding traffic modeling, the different way the road boundaries are treated makes the choice of the GWD lean towards PR in the phase of calibration and validation, i.e.\ when one has to estimate, at best, the model parameters, while it makes the choice lean toward FG when the model is used for traffic forecast, both for computational and modeling reasons: indeed, FG is the fastest in terms of CPU time and naturally focuses on the internal part of the road rather than on boundaries. This feature is useful since traffic forecast at boundaries is necessarily inaccurate due to the impossibility to estimate inflow and outflow traffic with good precision.

In the future, one could exploring the possibility to compute the considered GWDs in more complicated domain like $\R^2$, $\R^3$, and (road) networks.

\section*{Funding}
M.B.\ and E.C.\ would like to thank the Italian Ministry of University and Research (MUR) to support this research with funds coming from PRIN Project 2022 PNRR entitled ``Heterogeneity on the road - Modeling, analysis, control'', No. 2022XJ9SX.

E.C.\ would like to thank the Italian Ministry of University and Research (MUR) to support this research with funds coming from PRIN Project 2022 entitled ``Optimal control problems: analysis, approximation and applications'', No. 2022238YY5.

G.F.\ would like to thank the Italian Ministry of University and Research (MUR) to support this research with funds coming from PRIN Project 2022 entitled ``Variational analysis of complex systems in materials science, physics and biology'', No. 2022HKBF5C, CUP B53D23009290006.

This study was carried out within the Spoke 7 of the MOST -- Sustainable Mobility National Research Center and received funding from the European Union Next-Generation EU (PIANO NAZIONALE DI RIPRESA E RESILIENZA (PNRR) – MISSIONE 4 COMPONENTE 2, INVESTIMENTO 1.4 – D.D. 1033 17/06/2022, CN00000023). 

This manuscript reflects only the authors' views and opinions. Neither the European Union nor the European Commission can be considered responsible for them. 

E.C.\ and F.L.I.\ are funded by INdAM--GNCS Project, CUP E53C23001670001, entitled ``Numerical modeling and high-performance computing approaches for multiscale models of complex systems''.

F.L.I.\ is funded by Sapienza -- University of Rome, project ``Partial differential equations towards control theory and climate modeling'' (project code: RM124190DEC62D0A).

M.B., E.C., and F.L.I.\ are members of the e INdAM research group GNCS, while G.F.\ is member of the e INdAM research group GNAMPA.

\section*{Acknowledgments}
Authors want to thank Benedetto Piccoli, Francesco Rossi, Giuseppe Savar\'e, and Giacomo Enrico Sodini for the useful discussions, and Anna Th\"unen for the contribution given to Section \ref{sec:GWD-figalli-numerica}.

\section*{Use of AI tools declaration}
The authors declare they have not used Artificial Intelligence (AI) tools in the creation of this article.


\section*{Conflict of interest}
The authors declare there is no conflict of interest.

\appendix
\section{Appendices}\label{appendice}
\subsection{Numerical tests for the Figalli \& Gigli approach}\label{sec:GWD-figalli-numericaltest}

\subsubsection*{Test 1}
\begin{figure}[h!]
	\centering
	\includegraphics[width=0.49\textwidth]{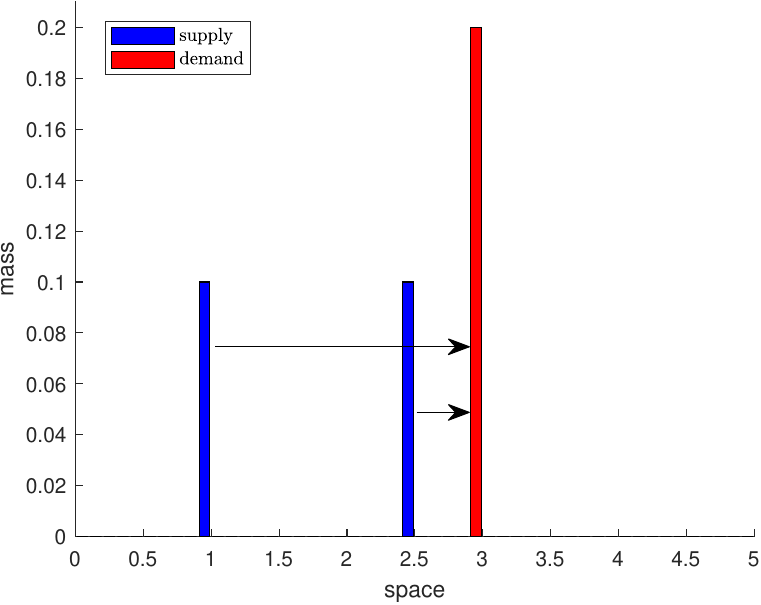}
	\includegraphics[width=0.49\textwidth]{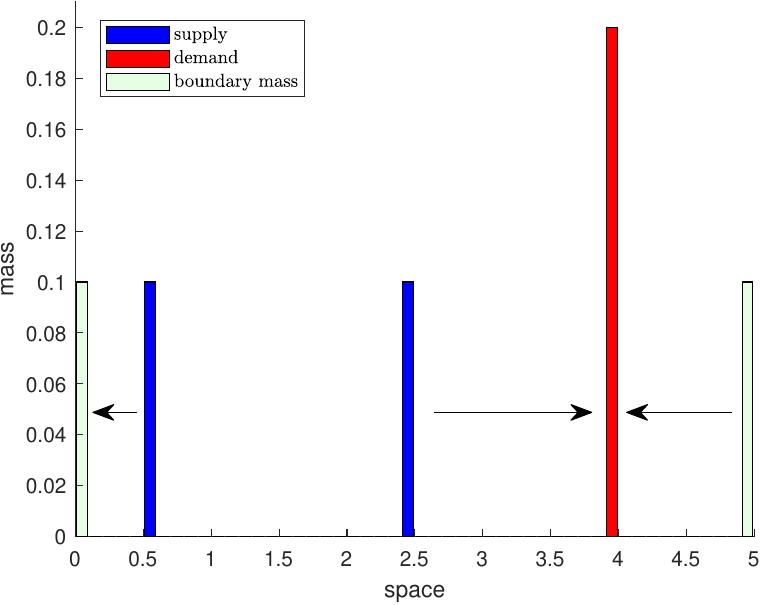}
	\caption{
		Test 1. Case 1.1 where boundaries are not used (left), and Case 1.2 where mass is both taken from and brought to boundaries (right). 
	}
	\label{fig:test-FG1}
\end{figure}

In this test we consider the case of atomic balanced masses in the domain $X=[0,5]$. We set $N=50$ so that $\Dx=0.1$. 
Total supply mass is 0.2, divided in two equal parts and concentrated in two points.
Total demand mass is also 0.2, concentrated in one point only.
\begin{enumerate}
	\item[1.1.] Here supply mass is located at $x=1$ and $x=2.5$, while demand mass is located at $x=3$, see Figure \ref{fig:test-FG1}(left). 
	Optimal mass transport map prescribes to move supply mass onto demand mass, without resorting to mass at boundaries. We have $\tilde W_1^\fg=0.25$.
	
	\item[1.2.] Here supply mass is located at $x=0.5$ and $x=2.5$, while demand mass is located at $x=4$, see Figure \ref{fig:test-FG1}(right). 
	Optimal mass transport map prescribes to move leftmost supply mass to the left boundary, rightmost supply mass onto the demand mass, and taking mass from right boundaries to bring them onto the demand mass. We have $\tilde W_1^\fg=0.3$.
\end{enumerate}

\subsubsection*{Test 2}
\begin{figure}[h!]
	\centering
	\includegraphics[width=0.5\textwidth]{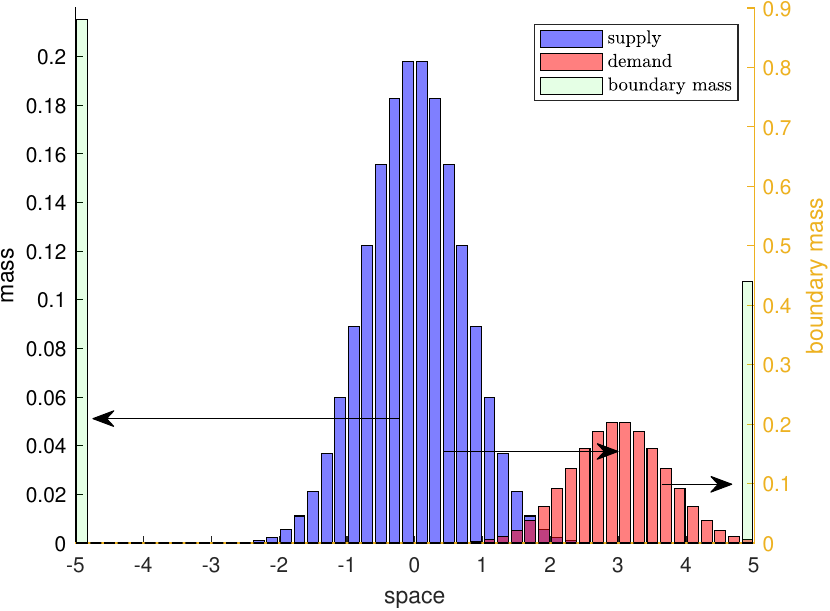}
	\caption{
		Test 2. Some supply (resp., demand) mass is brought to left (resp., right) boundary, the rest is transported.
	}
	\label{fig:test-FG2}
\end{figure}

In this test we consider the case of continuous unbalanced distributions in the domain $X=[-5,5]$,
$$
\rho^\supply(x) =e^{-x^2},\qquad
\rho^\demand(x) =\frac{e^{-(x-3)^2}}{4},
$$
see Figure \ref{fig:test-FG2}. We set $N=50$ so that $\Dx=0.2$. 

In this case the optimal mass transport map prescribes to move a part of the supply mass to the left boundary and a part of the demand mass to the right boundary, then transport the remaining mass. We have $\tilde W_1^\fg=6.84$.

\subsection{Numerical tests for the Piccoli \& Rossi approach}\label{sec:GWD-piccoli-numericaltest}

\subsubsection*{Test 1}
Consider the simple example proposed in \cite[Section 2.2]{piccoli2014ARMA}, where authors consider the two constant distributions 
$$
\rho^\supply(x) =\chi_{[-1,0]}(x),\qquad
\rho^\demand(x) = \chi_{[\xi,1+\xi]}(x),\qquad x\in \R
$$
with parameter $\xi\geq 0$.
The exact value of the GWD (with $a=b=p=1$) is easily computed and gives 
$$
W^\pr_1(\rho^\supply,\rho^\demand;1,1)=
\left\{
\begin{array}{ll}
	1+\xi-\xi^2/4, & \text{if } \xi\in[0,2], \\
	2, & \text{if }\xi\geq 2.
\end{array}
\right. 
$$

\begin{figure}[h!]
	\centering
	\includegraphics[width=0.5\textwidth]{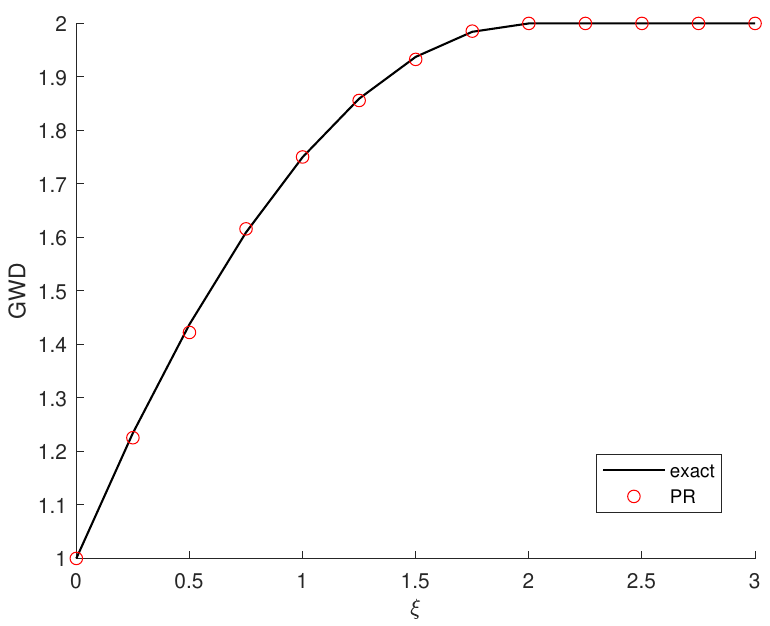}
	\caption{
		Test 1. Exact solution compared with numerical solution ($N=200$).
	}
	\label{fig:test-PR1}
\end{figure}

Figure \ref{fig:test-PR1} shows the result of the numerical computation in the domain $X=[-2,5]$ divided in $N=200$ cells, for $\xi\in[0,3]$. 

\subsubsection*{Test 2}
In this test we run the algorithm in a nontrivial case where exact solution is not known. More in detail, we consider $X=[-4,4]$, $N=100$, $a=b=p=1$, and
$$
\rho^\supply(x) =\frac{e^{-x+1}}{5}\chi_{[-2,0]}(x),\qquad
\rho^\demand(x) = e^{-(x-1)^2}.
$$
\begin{figure}[h!]
	\centering
	\includegraphics[width=0.49\textwidth]{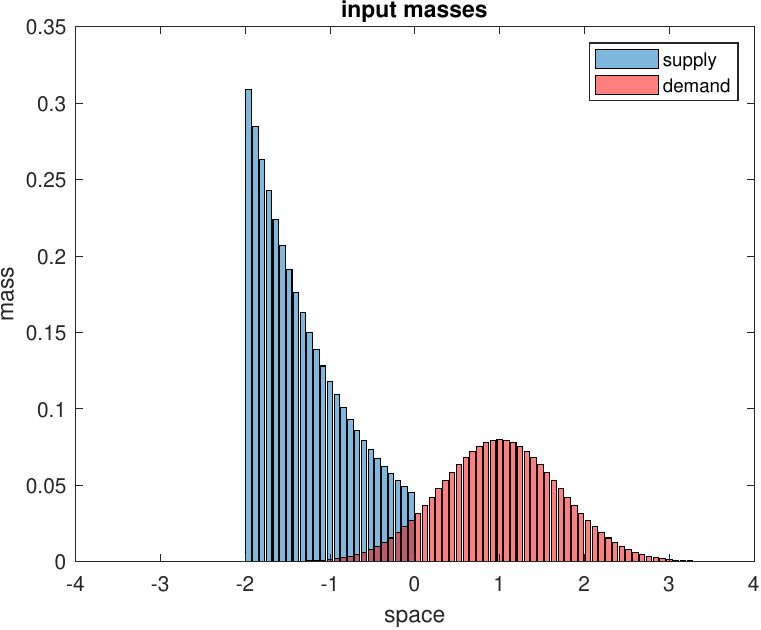}
	\includegraphics[width=0.49\textwidth]{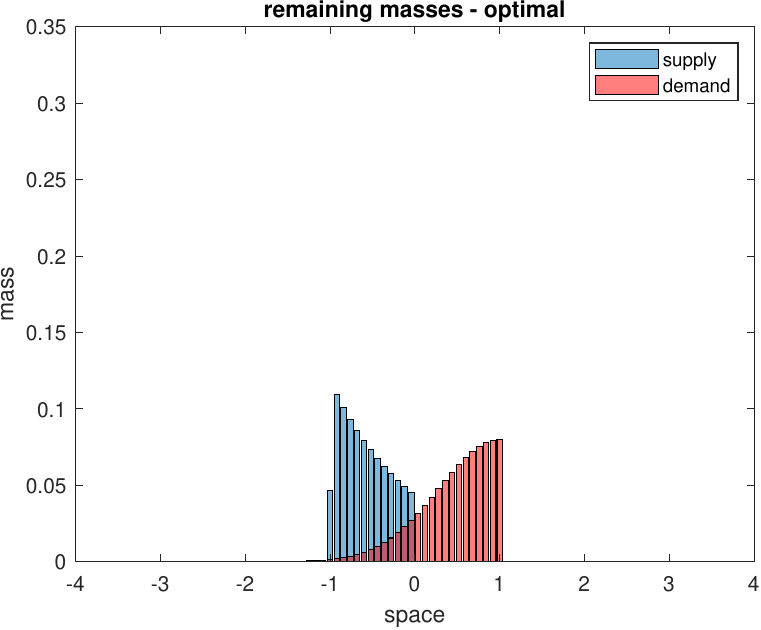}
	\caption{
		Test 2.	Initial distributions $m^\supply$ and $m^\demand$ (left) and remaining masses $\mathring m^\supply$ and $\mathring m^\demand$ corresponding to the optimal solution (right).
	}
	\label{fig:test-PR2}
\end{figure}
Figure \ref{fig:test-PR2} shows the result in terms of remaining mass after minimization. 
One can see that here the optimal solution ($\tilde W^\pr_1=4.36$) requires to destroy only a part of the mass.

\subsection{Numerical tests for the Gaussian Hellinger--Kantorovich and the Savaré \& Sodini approaches}\label{sec:GWD-SS-numericaltest}

We present here numerical tests for both the approaches described in Sects.\ \ref{sec:GHK-numerica} and \ref{sec:GWD-savare-numerica}. 
We recall that the two approaches, although quite different from the computational point of view, are theoretically equivalent (with $a=b=1$, which will be always the case from now on). 
This is also confirmed by our numerical experiments.

\subsubsection*{Test 1}
This test is mainly conceived to check the correctness of the code on simple examples with known exact solution. 
Let us consider the simplest possible scenario, with only one concentrated supply mass $\delta^\supply$ and only one concentrated demand mass $\delta^\demand$.
To simplify the notations, we denote by $m^\supply$, $x^\supply$ and by $m^\demand$, $x^\demand$ the mass and the position of the supply mass and the demand mass, respectively.
In this case the exact solution is 
$\tilde W^\ghk_2 = 
\tilde W^\sss_2 = 
H_g(x^\supply,m^\supply; x^\demand,m^\demand)$, 
where $H_g$ is defined in \eqref{def:H-SS}.

\begin{figure}[h!]
	\centering
	\includegraphics[width=0.49\textwidth]{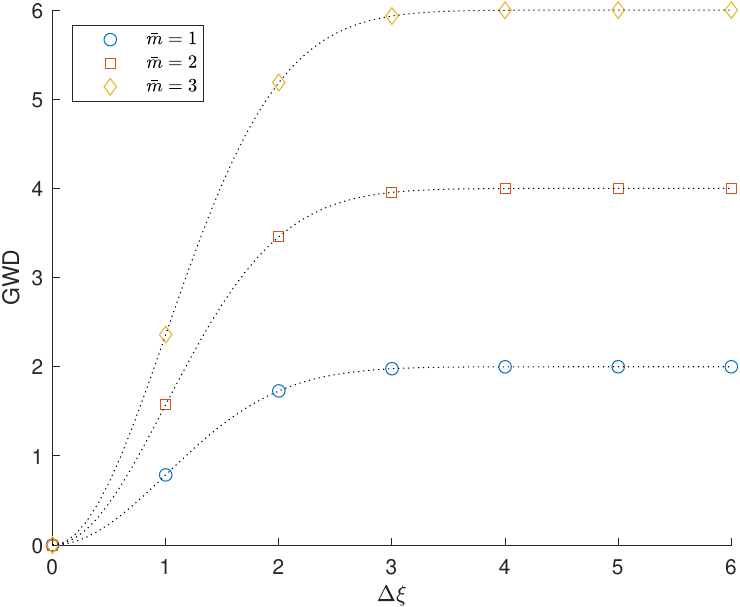}
	\includegraphics[width=0.49\textwidth]{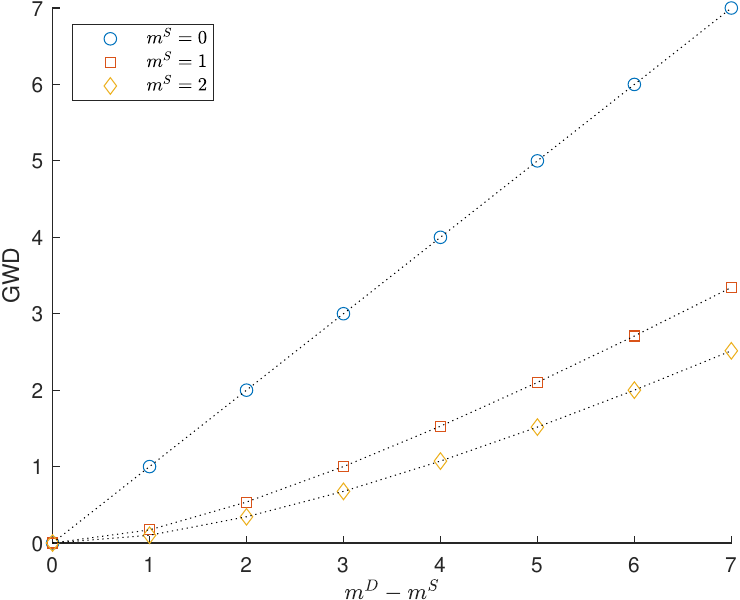}
	\caption{Test 1. GWD between two equal atomic mass distributions with mass $\bar m$ as a function of their distance $\Delta \xi$ (left), and GWD between two overlapping atomic mass distributions as a function of the difference of mass (right). 
		Dotted black line is the exact solution. 
	}
	\label{fig:test-SS1}
\end{figure}

Figure \ref{fig:test-SS1} shows the results of the algorithm for $\tilde W^\ghk_2$ as well as both algorithms for $\tilde W^\sss_2$ (SS-A and SS-B), alongside the exact solution.

More in detail, Figure \ref{fig:test-SS1}(left) shows the GWD in the case the two distributions are balanced ($m^\supply=m^\demand=\bar m$) and their distance progressively increases from 0 to 6. One can see that the GWD rapidly saturates at $2\bar m$.
As for $\tilde W^{\textsc{ghk}}_2$, that is consistent
with the fact that increasing the distance results
in relaxing the constraint \eqref{defcperGHK}
for the maximization of 
the dual functional, which involves
exponential functions:
from this, the maximum value
$2\bar m$ in the limit case of an infinite distance and
the exponential decay of the difference between this value and the actual maximum for finite distances that
grow large; as
for $\tilde W^{\textsc{ss}}_2$,
recalling \eqref{def:H-SS},
the objective in \eqref{def:S} decays
exponentially to $r_{jh}^\textsc{s}+r_{k\ell}^\textsc{d}$ which, in turn, equals at most
$2\bar m$ due to the constraints
on the homogeneous marginals.

Figure \ref{fig:test-SS1}(right) shows instead the GWD in the case the two distributions are located at the same point and the difference of their masses progressively increases from 0 to 7.

\subsubsection*{Test 2}
The goal of this test is mainly comparing the algorithms SS-A and SS-B.
We consider another simple scenario with only two concentrated masses for supply, and only one concentrated mass for demand, see Figure \ref{fig:test-SS2-scenario}.
\begin{figure}[h!]
	\centering
	\begin{tikzpicture}[scale=1.0]
	\draw[thick] (-8.5,0) -- (-0.5,0);
	\node at (-8.1,0) [above] {\tiny supply};
	\node at (-8.0,0) [below] {$x=$};
	\draw[thick] (0.5,0) -- (8.5,0);
	\node at (8.0,0) [above] {\tiny demand};
	\node at (1,0) [below] {$x=$};
	\node at (-5.75,0)[below]{0}; 
	\draw (-6,0) rectangle (-5.5,0.5); 
	\draw (-6,0.5) rectangle (-5.5,1); 
	\node at (-5.75,0.25){\footnotesize $r^\supply_{\cdot 1}$};
	\node at (-5.75,0.75){\footnotesize $r^\supply_{\cdot 2}$};
	\node at (-4.75,0)[below]{1}; 
	\draw (-5,0) rectangle (-4.5,0.5); 
	\node at (-4.75,0.25){\footnotesize $r^\supply_{\cdot 1}$};
	\node at (5.75,0)[below]{3}; 
	\draw (6,0) rectangle (5.5,0.5); 
	\draw (6,0.5) rectangle (5.5,1); 
	\draw (6,1) rectangle (5.5,1.5);
	\draw (6,1.5) rectangle (5.5,2);
	\node at (5.75,0.25){\footnotesize $r^\demand_{\cdot 1}$};
	\node at (5.75,0.75){\footnotesize $r^\demand_{\cdot 2}$};
	\node at (5.75,1.25){\footnotesize $r^\demand_{\cdot 3}$};
	\node at (5.75,1.75){\footnotesize $r^\demand_{\cdot 4}$};
\end{tikzpicture}
	\caption{Test 2: two concentrated masses for supply mass and one for demand mass.}
	\label{fig:test-SS2-scenario}
\end{figure} 
We assume that supply masses are located at $x=0$ and $x=1$, have masses $m=2$ and $m=1$, and have $n=2$, $n=1$ additional divisions, respectively. 
Demand mass is instead located at $x=3$, has mass $m=4$, and has $n=4$ additional divisions.
In this scenario it results $d=12$.  

\begin{figure}[h!]
	\centering
	\includegraphics[width=0.9\textwidth]{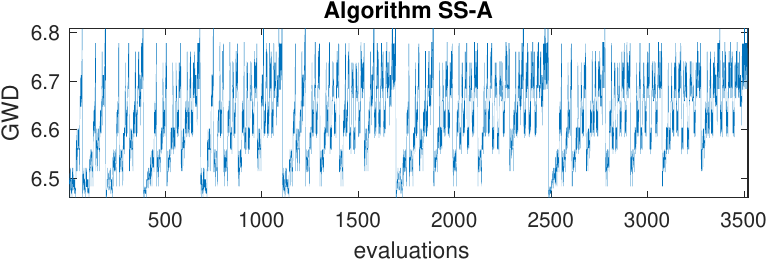}\\
	\includegraphics[width=0.9\textwidth]{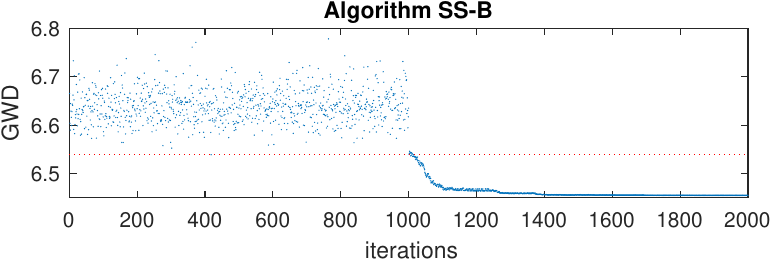}
	\caption{
		Test 2.	Algorithm SS-A, complete view of all the evaluations of $\tilde W^\sss_2$ (top), and 
		Algorithm SS-B, 1,000 random searches followed by about 1,000 descend iterations (bottom). Dotted red line indicates the minimum value found in the first phase of SS-B.
	}
	\label{fig:test-SS2}
\end{figure}

Figure \ref{fig:test-SS2}(top) shows the history of all the evaluations of $\tilde W^\sss_2$ done by the exhaustive algorithm SS-A, in the same order as they are computed, with $Q=6$. 
One can observe a recursive pattern due to the way $\Gamma$'s are generated (Step 2).
Minimum is already found among the first evaluations and it is $\tilde W^\sss_2=6.45$.

Figure \ref{fig:test-SS2}(bottom) shows the tentative values of $\tilde W^\sss_2$ found in the first exploratory phase of the algorithm SS-B (Steps 1-3), followed by the minimum values found during the random descend, with $Q=1,\!000$. 
The algorithm returns again $\tilde W^\sss_2=6.45$. 
Convergence seems to be well guaranteed and the two algorithms give the same value for $\tilde W^\sss_2$ ($\tilde W^\ghk_2$ is also the same).

For these reasons, in the following we will use algorithm SS-B only, which is obviously faster.

\bibliographystyle{plain}
\bibliography{biblio}

\begin{thebibliography}{10}

\bibitem{Amb-Mai-Ser2011}
L.~Ambrosio, Mainini E., and Serfaty S.
\newblock Gradient flow of the chapman-rubinstein-schatzman model for signed
  vortices.
\newblock {\em Ann. Inst. H. Poincar\'e C Anal. Non Lin\'eaire}, 28(2), 2011.

\bibitem{aw2000SIAP}
A.~Aw and M.~Rascle.
\newblock Resurrection of ``second order'' models of traffic flow?
\newblock {\em SIAM J. Appl. Math.}, 60:916--938, 2000.

\bibitem{BenBre2000NumerMath}
J.-D. Benamou and Y.~Brenier.
\newblock A computational fluid mechanics solutions to the
  {M}onge--{K}antoriovich mass transfer problem.
\newblock {\em Numer.\ Math.}, 84(3):375--393, 2000.

\bibitem{briani2018CMS}
M.~Briani, E.~Cristiani, and E.~Iacomini.
\newblock Sensitivity analysis of the {LWR} model for traffic forecast on large
  networks using {W}asserstein distance.
\newblock {\em Commun. Math. Sci.}, 16(1):123--144, 2018.

\bibitem{brianiAVpp}
M.~Briani, E.~Cristiani, and E.~Onofri.
\newblock Inverting the fundamental diagram and forecasting boundary
  conditions: {H}ow {M}achine {L}earning can improve macroscopic models for
  traffic flow.
\newblock To appear on Advances in Computational Mathematics. Preprint
  arXiv:2303.12740.

\bibitem{CafMcC2010ANNALS}
L.~Caffarelli and R.~J. {McCann}.
\newblock Free boundaries in optimal transport and {M}onge-{A}mp\`ere obstacle
  problems.
\newblock {\em Annals of Mathematics}, 171(2):673--730, 2010.

\bibitem{ChiPeySchVia2018FoundCompMath}
L.~Chizat, G.~Peyr\'e, B.~Schmitzer, and F.-X. Vialard.
\newblock An interpolating distance between optimal transport and
  {F}isher-{R}ao metrics.
\newblock {\em Foundation of Computational Mathematics}, 18:1--44, 2018.

\bibitem{ChiPeySchVia2018JFA}
L.~Chizat, G.~Peyr\'e, B.~Schmitzer, and F.-X. Vialard.
\newblock Unbalanced optimal transport: {D}ynamic and {K}antorovich
  formulations.
\newblock {\em Journal of Functional Analysis}, 274:3090--3123, 2018.

\bibitem{cristianibook}
E.~Cristiani, B.~Piccoli, and A.~Tosin.
\newblock {\em Multiscale modeling of pedestrian dynamics}, volume~12 of {\em
  Series Modeling, Simulation \& Applications}.
\newblock Springer, 2014.

\bibitem{cristiani2019M2AS}
E.~Cristiani and M.~C. Saladino.
\newblock Comparing comparisons between vehicular traffic states in microscopic
  and macroscopic first-order models.
\newblock {\em Mathematical Methods in the Applied Sciences}, 42(3):918--934,
  2019.

\bibitem{Figalli2010ARMA}
A.~Figalli.
\newblock The optimal partial transport problem.
\newblock {\em Archive for Rational Mechanics and Analysis}, 195(2):533--560,
  2010.

\bibitem{figalli2010JMPA}
A.~Figalli and N.~Gigli.
\newblock A new transportation distance between non-negative measures, with
  applications to gradients flows with {D}irichlet boundary conditions.
\newblock {\em Journal de Math{\'e}matiques Pures et Appliqu{\'e}es},
  94(2):107--130, 2010.

\bibitem{Hanin1992PAMS}
L.~G. Hanin.
\newblock {K}antorovich-{R}ubinstein norm and its application in the theory of
  {L}ipschitz spaces.
\newblock {\em Proceedings of the American Mathematical Society},
  115(2):345--352, 1992.

\bibitem{KantorovichRubinstein1958VestnLeningrUniv}
L.~Kantorovich and G.~S. Rubinstein.
\newblock On a space of totally additive functions.
\newblock {\em Vestnik Leningradskogo Universiteta}, 13:52--59, 1958.

\bibitem{KonMorVor2016AdvDE}
S.~Kondratyev, L.~Monsaingeon, and D.~Vortnikov.
\newblock A new optimal transport distance on the space of finite {R}adon
  measures.
\newblock {\em Adv.\ Differential Equations}, 21(11/12):1117--1164, 2016.

\bibitem{LieMieSav2018InventMath}
M.~Liero, A.~Mielke, and G.~Savar\'e.
\newblock Optimal entropy-transport problems and a new
  {H}ellinger-{K}antorovich distance between positive measures.
\newblock {\em Invent.\ Math.}, 211(3):969--1117, 2018.

\bibitem{lighthill1955PRSLA}
M.~J. Lighthill and G.~B. Whitham.
\newblock On kinematic waves {II}. {A} theory of traffic flow on long crowded
  roads.
\newblock {\em Proc. R. Soc. Lond. Ser. A}, 229:317--345, 1955.

\bibitem{lombardini2022PAMS}
L.~Lombardini and F.~Rossi.
\newblock Obstructions to extension of {W}asserstein distances for variable
  masses.
\newblock {\em Proceedings of the American Mathematical Society},
  150(11):4879--4890, 2022.

\bibitem{ma2021AAAI}
Z.~Ma, X.~Wei, X.~Hong, H.~Lin, Y.~Qiu, and Y.~Gong.
\newblock Learning to count via unbalanced optimal transport.
\newblock In {\em Proceedings of the AAAI Conference on Artificial
  Intelligence}, volume~35, pages 2319--2327, 2021.

\bibitem{Mai2012}
E.~Mainini.
\newblock A description of transport cost for signed measures.
\newblock {\em J. Math. Sci. (N.Y.)}, 181(6), 2012.

\bibitem{piccoli2014ARMA}
B.~Piccoli and F.~Rossi.
\newblock Generalized {W}asserstein distance and its application to transport
  equations with source.
\newblock {\em Archive for Rational Mechanics and Analysis}, 211:335--358,
  2014.

\bibitem{piccoli2016ARMA}
B.~Piccoli and F.~Rossi.
\newblock On properties of the generalized {W}asserstein distance.
\newblock {\em Arch. Rational Mech. Anal.}, 222:1339--1365, 2016.

\bibitem{richards1956OR}
P.~I. Richards.
\newblock Shock waves on the highway.
\newblock {\em Oper. Res.}, 4:42--51, 1956.

\bibitem{santambrogiobook}
F.~Santambrogio.
\newblock {\em Optimal Transport for Applied Mathematicians. Calculus of
  Variations, PDEs, and Modeling}, volume~87 of {\em Series Progress in
  Nonlinear Differential Equations and Their Applications}.
\newblock Birkh\"auser, 2015.

\bibitem{savare2024JMPA}
G.~Savar{\'e} and G.~E. Sodini.
\newblock A relaxation viewpoint to {U}nbalanced {O}ptimal {T}ransport:
  {D}uality, optimality and {M}onge formulation.
\newblock {\em Journal de Math{\'e}matiques Pures et Appliqu{\'e}es},
  188:114--178, 2024.

\bibitem{sinhabook}
S.~M. Sinha.
\newblock {\em Mathematical Programming. Theory and Methods}.
\newblock Elsevier, 2006.

\bibitem{villani2003book}
C.~Villani.
\newblock {\em Topics in Optimal Transportation}, volume~58 of {\em Series
  Graduate Studies in Mathematics}.
\newblock American Mathematical Society, 2003.

\bibitem{villani2009book}
C.~Villani.
\newblock {\em Optimal Transport. Old and New}, volume 338 of {\em Series
  Grundlehren der mathematischen Wissenschaften}.
\newblock Springer, 2009.

\bibitem{zhang2002TRB}
H.~M. Zhang.
\newblock A non-equilibrium traffic model devoid of gas-like behavior.
\newblock {\em Transportation Res. Part B}, 36:275--290, 2002.

\end{thebibliography}


\end{document}